\documentclass[a4paper,11pt]{article}

\usepackage{authblk}
\usepackage{parskip}
\usepackage{mathtools}
\usepackage{amssymb}
\usepackage{amsmath}
\usepackage{latexsym}
\usepackage{mathrsfs}  
\usepackage{bbm}       
\usepackage{amsthm}    
\usepackage{relsize}   
\usepackage{graphicx}
\usepackage{thmtools}  
\usepackage{mathrsfs}  
\usepackage{paralist}  
\usepackage{lipsum}    
\usepackage[all]{xy}   

\numberwithin{equation}{section}

\usepackage{titlesec}   

\titleformat{\section}[block]{\bfseries\filcenter}
{{\upshape\thesection\enspace}}{.5em}{}

\titleformat{\subsection}[block]{\filcenter}
{{\upshape\thesubsection\enspace}}{.5em}{} 

\usepackage{enumitem}  
\setlist{nosep}  
\setitemize[0]{leftmargin=*}
\setenumerate[0]{leftmargin=*}
\setenumerate[1]{label={(\arabic*)}} 

\newcommand{\N}{\mathbb{N}}     
\newcommand{\R}{\mathbb{R}}     
\newcommand{\Prob}{\mathbb{P}}  
\newcommand{\Exp}{\mathbb{E}}   
\newcommand*{\backin}{\rotatebox[origin=c]{-180}{ $\in$ }}%
\newcommand{\st}{\,:\,}         
\newcommand{\goth}[1]{\mathfrak{#1}} 
\newcommand{\ind}[2]{\mathbbm{1}_{#1}\left( #2 \right)}          
\newcommand{\inner}[2]{\left\langle #1 \, , \, #2 \right\rangle} 
\newcommand{\norm}[1]{\left|\left|#1\right|\right|}              
\newcommand{\triplet}[3]{\left( #1, #2, #3 \right) }             
\newcommand{\ProbSpace}{\triplet{\Omega}{\mathscr{F}}{\Prob}}    
\newcommand{\abs}[1]{\left| #1 \right|}                          
\renewcommand{\qedsymbol}{$\square$}                       

\newcommand{\defeq}{\mathrel{\mathop:}=}                         
\newcommand\restr[2]{{
  \left.\kern-\nulldelimiterspace 
  #1 
  \vphantom{\big|} 
  \right|_{#2} 
  }}
  
\theoremstyle{plain} 
\newtheorem{theo}{Theorem}[section]    
\newtheorem{prop}[theo]{Proposition} 
\newtheorem{coro}[theo]{Corollary}
\newtheorem{lemm}[theo]{Lemma}
\newtheorem{assu}[theo]{Assumption}
\newtheorem{rema}[theo]{Remark}

\theoremstyle{definition} 
\newtheorem{defi}[theo]{Definition}
\newtheorem{exam}[theo]{Example}

\newtheorem{nota}[theo]{Notation}

\declaretheoremstyle[%
  spaceabove=-5pt,%
  spacebelow=6pt,%
  headfont=\normalfont\itshape,%
  postheadspace=1em,%
  qed=\qedsymbol%
]{mystyle} 
\declaretheorem[name={Proof},style=mystyle,unnumbered,
]{prf}

 \title{Stochastic integration in Hilbert spaces with respect to cylindrical martingale-valued measures}

\author{
A. E. Alvarado-Solano\thanks{anddy.alvarado@ucr.ac.cr }   }

\author{   
C. A. Fonseca-Mora\thanks{christianandres.fonseca@ucr.ac.cr}} 
\affil{Escuela de Matem\'{a}tica, Universidad de Costa Rica, 
\\ San Jos\'{e}, 11501-2060, Costa Rica}
\date{}   
\begin{document}

\maketitle

\abstract{In this work we introduce a theory of stochastic integration for operator-valued integrands with respect to some classes of cylindrical martingale-valued measures in Hilbert spaces.  The integral is constructed via the radonification of cylindrical martingales by a Hilbert-Schmidt operator theorem and unifies several other theories of stochastic integration in Hilbert spaces. In particular, our theory covers the theory of stochastic integration with respect to a Hilbert space valued L\'{e}vy process (which is not required to satisfy any moment condition), with respect to a cylindrical L\'{e}vy processes with (weak) second moments and with respect to a L\'{e}vy-valued random martingale measures with finite second moment. As an application of our theory of integration we prove existence and uniqueness of solutions for stochastic stochastic partial differential equations  driven by multiplicative cylindrical martingale-valued measure noise with rather general coefficients.}

\smallskip



\section{Introduction}\label{sectIntro}

In recent years, there has been an increasing interest into the study of stochastic partial differential equations driven by cylindrical noise (e.g. \cite{BrzezniakZabczyk:2010, KosmalaRiedle:VariaSPDE, KumarRiedle:StochCaucLevy, LiuZhai:2016, PriolaZabczyk, Riedle:2015, Riedle:2018}). Motivated by these developments, in this work we introduce a construction for the stochastic integral with respect to some classes of cylindrical martingale-valued measures. Our theory includes the cases where the integrator is a Hilbert space valued L\'{e}vy processes (with no moments assumptions), a cylindrical L\'{e}vy process with (weak) finite second moments in a Hilbert space, or a L\'{e}vy-valued random martingale measure with finite second moment. 

The concept of cylindrical martingale-valued measures in a locally convex space was introduced in the work \cite{FonsecaMora:2018-1}, and generalizes to locally convex spaces the martingale-valued measures
introduced for the finite dimensional setting by Walsh \cite{Walsh:1986} and then extended to infinite dimensional settings such as Hilbert spaces \cite{Applebaum:2006} and duals of nuclear Fréchet spaces \cite{Xie:2001}. However, in this work we will consider such objects only in the Hilbert space setting. 

Roughly speaking, a cylindrical martingale-valued measure in a Hilbert space $H$ is a family $M=(M(t,A): t \geq 0, A \in \mathcal{R})$ such that $(M(t,A): t \geq 0)$ is a cylindrical martingale in $H$ for each $A \in \mathcal{R}$ and $M(t,\cdot)$ is finitely additive on $\mathcal{R}$ for each $t \geq 0$. Here $\mathcal{R}$ is a ring of Borel subsets of a topological space $U$. We will also request the existence of some weak second moments for $M$, which are determined by a family of continuous Hilbertian semi-norms $(q_{r,u}: r \geq 0, u \in U)$ on $H$ (see Definition \ref{defiMartingValuedMeasure}). Typical examples include (but are not limited to) the martingale part of a $H$-valued L\'{e}vy process (see Example \ref{examLevyMartValuedMeasure}), the cylindrical L\'{e}vy processes which possesses (weak) second moments (see Example \ref{examCylLevySquareMarting}), and the L\'{e}vy-valued random martingale measures (see Example \ref{examLevyValRandMartMeasu}) of which Gaussian space-time white noise is a particular case. 

Our main task will be to develop a theory of stochastic integration for operator-valued integrands with respect to a given cylindrical martingale-valued measure $M$. Our integrands are families of random operators $(\Phi(r,\omega,u): r \geq 0, \omega \in \Omega, u \in U)$ such that each $\Phi(r,\omega, u)$ is Hilbert-Schmidt from a Hilbert space $H_{q_{r,u}}$ determined by $H$ and the continuous Hilbertian seminorm $q_{r,u}$ on $H$ into a Hilbert space $G$, and such that this family satisfies a square integrability condition (see Definition \ref{defiIntegrandsWithSquareMoments}). In the case of integration with respect to a $H$-valued L\'{e}vy process our square integrability condition coincides with those obtained for example in \cite{Applebaum:2006, PeszatZabczyk} (see Example \ref{examStochIntegForLevyProcess}), and in the case of integration with respect to a cylindrical L\'{e}vy process in $H$ with (weak) second moments our square integrability  condition allows for a larger class of integrands than those  in \cite{Riedle:2014} and coincides with that in \cite{KosmalaRiedle:VariaSPDE} (see  Example \ref{examItoIsomeCylLevySquaInteg}). The relationship of these theories of stochastic integration with the stochastic integral developed in this work will be explored thoroughly (see Section \ref{sectExamStochInte}). 

%
%

The construction of the stochastic integral follows an It\^{o}'s approach by first defining it on a
class of simple processes and then extending it via an It\^{o} isometry. In the case of simple processes, the construction of the stochastic integral is carried out via the radonification of cylindrical martingales by a Hilbert-Schmidt operator theorem. 

To explain our approach in more detail, suppose $H$ and $G$ are separable Hilbert spaces and let  $M=(M(t,A): t \geq 0, A \in \mathcal{R})$ be a $H$-valued martingale-valued measure defined on a probability space $\ProbSpace$ equipped with a filtration $(\mathcal{F}_{t})$. The simplest class of integrands is of the form $\Phi(r,\omega,u)= \mathbbm{1}_{]a,b]}(r) \mathbbm{1}_{F}(\omega) \mathbbm{1}_{A}(u) S$,  
where $0 \leq a< b \leq T$, $F \in \mathcal{F}_{a}$, $A \in \mathcal{R}$ and $S$ is a Hilbert-Schmidt operator from $H$ into $G$. Any sensible definition of the stochastic integral leads to the equality:
\begin{eqnarray*}
\inner{\int_{0}^{t}\int_{U}  \Phi(r,u) M(dr,du)}{g} 
& = & \inner{\mathbbm{1}_{F}S(M(b \wedge t,A)-M(a \wedge t,A))}{g} \\
& = & \mathbbm{1}_{F} \inner{M(b \wedge t,A)-M(a \wedge t,A)}{S^{*}g},
\end{eqnarray*}
for every $g \in G$ and $t \in [0,T]$, and where $S^{*}$ denotes the adjoint operator of $S$. If $M$ is assumed to be a cylindrical martingale-valued measure with (weak) second moments in $H$, then for each $r \geq 0$ and $A \in \mathcal{R}$ we have by definition that $M(r,A)$  is a continuous linear operator from $H$ into the space $L^{2} \ProbSpace$ of square integrable real-valued random variables. Hence, in the equality above we must substitute the inner product by the application of the linear operator  $M(b \wedge t,A)-M(a \wedge t,A)$ to $S^{*}g $; then we arrive at:
\begin{equation} \label{eqIntroDescriInteg}
\inner{\int_{0}^{t}\int_{U} \Phi(r,u) M(dr,du)}{g} = \mathbbm{1}_{F} (M(b \wedge t,A)-M(a \wedge t,A))(S^{*}g),
\end{equation}
for every $g \in G$ and $t \in [0,T]$. If $S$ is only a continuous linear operator the above equality defines the integral as a cylindrical process in $H$. However, since we are assuming that $S$ is a Hilbert-Schmidt operator from $H$ into $G$ the radonification by a single Hilbert-Schmidt theorem (see Theorem \ref{theoRadonCylMartHilbSchmi}) proves that the equality above defines a $G$-valued square integrable martingale with (strong) second moments which we can define as our integral. Here is where our approach differs from that in  \cite{JakubowskiRiedle:2017, KosmalaRiedle:VariaSPDE, Riedle:2014} because in these works the authors consider only radonification of the increments, i.e. the radonification in \eqref{eqIntroDescriInteg} occurs for a fixed $t \in [0,T]$ and not for the whole process as we do in this work. 

Our assumption of (weak) second moments on the cylindrical martingale-valued measure $M$ permits that our simple integrands satisfy an It\^{o} isometry and hence the extension to integrands with second moments follows from standard arguments. We will also prove that our stochastic integral satisfies all the standard properties and comparisons with other theories of integration will be discussed.

 As an application of our theory of stochastic integration we consider existence and uniqueness of mild solutions to  stochastic partial differential equations driven by cylindrical martingale-valued measures under some Lipschitz and growth conditions. Existence and uniqueness of solutions in the presence of multiplicative L\'{e}vy noise (with no moments assumptions) is also proved.  

The organization of the paper is as follows. In Sect. \ref{sectionPrelim} we review the most relevant concepts and results on cylindrical processes in Hilbert spaces that we will need throughout this work. In particular, we prove the regularization by a single Hilbert-Schmidt operator theorem for cylindrical martingales with weak second moments. In Sect. \ref{SectCylMarValMeasu} we review the definition of cylindrical martingale-valued measures in Hilbert spaces and present several examples. The construction of the stochastic integral, the study of its properties, as well as the presentation of examples and comparisons with the literature is carried out in Sect. \ref{sectStochaIntegral}. Finally, in Sect. \ref{sectStocEvoluEquations} we study the existence and uniqueness of solutions to stochastic partial differential equations driven by cylindrical martingale-valued measures.

\section{Preliminaries} \label{sectionPrelim}

Let $H$ be a separable Hilbert space with inner product $\inner{\cdot}{ \cdot }$ and norm $\norm{\cdot}$. We denote by $\mathcal{B}(H)$ its Borel $\sigma$-algebra. We always identify $H$ with its (strong) dual space.  For a $\sigma$-finite measure space $(S, \mathcal{S}, \mu)$, denote by $L^{0}( S, \mathcal{S}, \mu; H)$ (or $L^{0}( S, \mathcal{S}, \mu)$ if $H=\R$) the linear space of all the equivalence classes of $H$-valued $\mathcal{F}/\mathcal{B}(H)$-measurable functions (or random variables if $\mu$ is a probability measure). We equip $L^{0}( S, \mathcal{S}, \mu; H)$ with the topology of convergence in $\mu$-measure and in this case it is a complete, metrizable, topological vector space. Similarly, for any $p \geq 1$ we denote by  $L^{p}(S, \mathcal{S}, \mu; H)$ the linear space of all the equivalence classes of $H$-valued measurable functions $X$ for which $\norm{X}_{L^{p}( S, \mathcal{S}, \mu; H)}^{p}= \int_{S} \, \norm{X}^{p} d \mu < \infty$ (in this case $X$ is called $p$-integrable).  This is a Banach space (Hilbert space if $p=2$) when equipped with the norm $\norm{\cdot}_{L^{p}( S, \mathcal{S}, \mu; H)}$. 

Let $G$ be another separable Hilbert space. The collection of all the continuous linear operators from $H$ into $G$ will be denoted by $\mathcal{L}(H,G)$. Recall that $S \in \mathcal{L}(H,G)$ is called Hilbert-Schmidt if for some (equivalently for any) complete orthonormal system $(h_{n}: n \in \N)$ in $H$, we have $\displaystyle{ \sum_{n=1}^{\infty} \norm{Sh_{n}}^{2}< \infty}$. The space of all the Hilbert-Schmidt operators is denoted by $\mathcal{L}_{2}(H,G)$ and it is a Hilbert space when equipped with the inner product and corresponding norm:
$$ \inner{R}{S}_{\mathcal{L}_{2}(H,G)}= \sum_{k=1}^{\infty} \inner{R h_{n}}{S h_{n}}, \quad  
\norm{S}_{\mathcal{L}_{2}(H,G)}= \left( \sum_{k=1}^{\infty} \norm{S h_{n}}^{2} \right)^{1/2}. $$

If $q$ is a continuous Hilbertian seminorm on $H$, we denote by $H_{q}$ the Hilbert space that corresponds to the completion of the pre-Hilbert space $(H / \mbox{ker}(q), \overline{q})$, where $\overline{q}(h+\mbox{ker}(q))=q(h)$ for each $h \in H$. The quotient  map $H \rightarrow H / \mbox{ker}(q)$ has a unique continuous linear extension $i_{q}: H \rightarrow H_{q}$. The Hilbertian norm in $H_{q}$ will be again denoted by $q$ and observe that $q(i_{q} h)=q(h)$ for each $h \in H$. 

Let $p$ be another continuous Hilbertian semi-norm on $H$ for which $q \leq p$. In this case, $\mbox{ker}(p) \subseteq \mbox{ker}(q)$ and the inclusion map from $H / \mbox{ker}(p)$ into $H / \mbox{ker}(q)$ is linear and continuous,  therefore it has a unique continuous extension $i_{q,p}:H_{p} \rightarrow H_{q}$. Moreover, we have $i_{q}=i_{q,p} \circ i_{p}$.

All our random variables will be considered defined on some given complete probability space $\ProbSpace$  equipped with a filtration $(\mathcal{F}_{t} : t \geq 0)$ that satisfies the \emph{usual conditions}, i.e. it is right continuous and $\mathcal{F}_{0}$ contains all subsets of sets of $\mathcal{F}$ of $\Prob$-measure zero. We denote by $\mathcal{P}_{\infty}$ the predictable $\sigma$-algebra on $[0, \infty) \times \Omega$ and for any $T>0$, we denote by $\mathcal{P}_{T}$ the restriction of $\mathcal{P}_{\infty}$ to $[0,T] \times \Omega$. 

Given $T>0$ we denote by $\mathcal{M}_{T}^{2}(H)$ the space of all the $H$-valued c\`adl\`ag square integrable martingales defined on $[0,T]$, which is a Banach space when equipped with the norm (see \cite{DaPratoZabczyk}, Proposition 3.9, p.79),
$$\norm{M}_{\mathcal{M}_{T}^{2}(H)}= \left( \Exp \sup_{0 \leq t \leq T} \norm{M_{t}}^{2} \right)^{1/2}. $$

For any $n \in \N$ and any $h_{1}, \dots, h_{n} \in H$, we define a linear map $\pi_{h_{1}, \dots, h_{n}}: H \rightarrow \R^{n}$ by $\pi_{h_{1}, \dots, h_{n}}(y)=(\inner{y}{h_{1}}, \dots, \inner{y}{h_{n}})$ $\forall \, y \in H$.  A cylindrical set in $H$ based on $M \subseteq H$ is a set of the form 
$$ \mathcal{Z}\left(h_{1}, \dots, h_{n}; A \right) = \left\{ y \in H \st (\inner{y}{h_{1}}, \dots, \inner{y}{h_{n}}) \in A \right\}= \pi_{h_{1}, \dots, h_{n}}^{-1}(A), $$
where $n \in \N$, $h_{1}, \dots, h_{n} \in M$ and $A \in \mathcal{B}\left(\R^{n}\right)$ The collection of all the cylindrical sets based on $M$ is denoted by $\mathcal{Z}(H,M)$. It is an algebra but if $M$ is a finite set then it is a $\sigma$-algebra. The $\sigma$-algebra generated by $\mathcal{Z}(H,M)$ is denoted by $\mathcal{C}(H,M)$. If $M=H$, we write $\mathcal{Z}(H)=\mathcal{Z}(H,M)$ and $\mathcal{C}(H)=\mathcal{C}(H,M)$. 
A function $\mu: \mathcal{Z}(H) \rightarrow [0,\infty]$ is called a \emph{cylindrical measure} on $H$, if for each finite subset $M \subseteq H$ the restriction of $\mu$ to $\mathcal{C}(H,M)$ is a measure. A cylindrical measure $\mu$ is said to be \emph{finite} if $\mu(H)< \infty$ and a \emph{cylindrical probability measure} if $\mu(H)=1$. 

A \emph{cylidrical random variable} in $H$ is a linear map $X: H \rightarrow L^{0} \ProbSpace$. We can associate to $X$ a cylindrical probability measure $\mu_{X}$ on $H$, called its \emph{cylindrical distribution}, and given by 
\begin{equation*} 
\mu_{X}(Z) \defeq \Prob \left( ( X(\phi_{1}), \dots, X(\phi_{n})) \in A  \right), 
\end{equation*}
for the cylindrical set $Z=\mathcal{Z}\left(h_{1}, \dots, h_{n}; A \right)$. 
Conversely, to every cylindrical probability measure $\mu$ on $H$ there is a canonical cylindrical random variable for which $\mu$ is its cylindrical distribution. 

Let $J=\R_{+} \defeq [0,\infty)$ or $J=[0,T]$ for some $T>0$. We say that $X=( X_{t}: t \in J)$ is a \emph{cylindrical process} in $H$ if $X_{t}$ is a cylindrical random variable for each $t \in J$. We say that $X$ is \emph{$p$-integrable} if $ \Exp \left( \abs{X_{t}(h)}^{p} \right)< \infty$, $\forall \, h \in H$ and $t \in J$, and has \emph{zero-mean} if $ \Exp \left( X_{t}(h) \right)=0$,  $\forall \, h \in H$ and $t \in J$. 

To every $H$-valued stochastic process $X=( X_{t} : t \in J)$ there corresponds a cylindrical process defined by $( \inner{X_{t}}{h} : t \in J)$, for each $h \in H$. We will say that it is the \emph{cylindrical process determined/induced} by $X$. Conversely, if a $H$-valued processes $Y=(Y_{t}: t \in J)$ is said to be a $H$-valued \emph{version} of the cylindrical process $X=(X_{t}: t \in J)$ if for each $t \in J$ and $h \in H$, $\inner{Y_{t}}{h}=X_{t}(h)$ $\Prob$-a.e.  

A cylindrical process $X=(X_{t}: t \in J)$ is called a \emph{cylindrical square integrable martingale}  if for each $h \in H$, $(X_{t}(h): t \in J)$ is a real-valued square integrable martingale. 

The following result concerns the construction of a Hilbert space valued square integrable martingale from a cylindrical square integrable martingale. This is often known as a \emph{radonification by a single Hilbert-Schmidt operator}, and it can be extended to the more general context of multi-Hilbertian spaces (see \cite{FonsecaMora:ReguLCS}). However, since this result will be of fundamental importance in our construction of the stochastic integral, for the benefit of the reader we decided to include a short self-contained proof of it. 

\begin{theo}\label{theoRadonCylMartHilbSchmi} Let $H$ and $G$ be separable Hilbert spaces and let $M=(M_{t}: t \in [0,T])$ be a cylindrical square integrable martingale in $H$ such that for each $t \in [0,T]$, the map $M_{t}: H \rightarrow L^{0} \ProbSpace$ is continuous. If $S \in \mathcal{L}_{2}(H,G)$, the cylindrical square integrable martingale  $M\circ S^{*}= ( M_{t}\circ S^{*} : t \in [0,T])$ has a $G$-valued c\`{a}dl\`{a}g version $X= (X_{t}: t \in [0,T])$ that is a square integrable martingale. 
\end{theo}
\begin{prf}
Let $\widetilde{M}=(\widetilde{M}_{t}: t \in [0,T])$ be such that for each $h \in H$, $(\widetilde{M}_{t}(h): t \in [0,T])$ is a c\`{a}dl\`{a}g version of $(M_{t}(h): t \in [0,T])$. It is easy to verify that $\widetilde{M}$ is a cylindrical square integrable martingale. Moreover for each $t \in [0,T]$, since $\widetilde{M}_{t}$ and $M_{t}$ share the same characteristic function, then the fact that $M_{t}: H \rightarrow L^{0} \ProbSpace$ is continuous implies that $\widetilde{M}_{t}: H \rightarrow L^{0} \ProbSpace$ is continuous. 

We claim that $\widetilde{M}: H \rightarrow \mathcal{M}^{2}_{T}(\R)$, defined by $h \mapsto \widetilde{M}(h)$, is linear and continuous. In effect, the mapping $\widetilde{M}$ is well-defined and obviously linear. To prove its continuity, by the closed graph theorem it is enough to show that it is closed. Let $h_{n} \rightarrow h$ in $H$, and suppose $\widetilde{M}(h_{n}) \rightarrow N$ in  $\mathcal{M}^{2}_{T}(\R)$. Then, for each $t \in [0,T]$, $\widetilde{M}_{t}(h_{n}) \rightarrow N_{t}$ in probability. But since $h_{n} \rightarrow h$, we also have $\widetilde{M}_{t}(h_{n}) \rightarrow \widetilde{M}_{t}(h)$ in probability. By uniqueness of limits $N_{t}= \widetilde{M}_{t}(h)$ $\Prob$-a.e. Since both processes are c\`{a}dl\`{a}g, then $N$ and $\widetilde{M}(h)$ are indistinguishable. Thus the mapping $\widetilde{M}$ is closed, hence continuous. Observe that the above also implies that for each $t \in [0,T]$, $\widetilde{M}_{t} \in \mathcal{L}(H, L^{2} \ProbSpace )$. 

Now since $S \in \mathcal{L}_{2}(H,G)$, there exist complete orthonormal systems $(h_{n}: n \in \N)$ in $H$, $(g_{n}: n \in \N)$ in $G$, and a sequence $(\lambda_{n}:n \in \N) \subseteq \R$ such that $\displaystyle{\norm{S}_{\mathcal{L}_{2}(H,G)}^{2}= \sum_{n=1}^{\infty} \abs{\lambda_{n}}^{2} < \infty}$, and 
$$ S h= \sum_{n=1}^{\infty} \lambda_{n} \inner{h_{n}}{h} g_{n}, \quad \forall \,  h \in H. $$
For each $n \in \N$, let 
$$  X^{(n)}= \sum_{k=1}^{n} \lambda_{k} \widetilde{M}(h_{k}) g_{k} \in \mathcal{M}^{2}_{T}(G). $$
We will show that the sequence $(X^{(n)}:n \in \N)$ is Cauchy in $\mathcal{M}^{2}_{T}(G)$. In effect, for  $m \geq n$, from Doob's inequality, and the arguments give in the above paragraphs we have 
\begin{eqnarray*}
\Exp \left( \sup_{0 \leq t \leq T} \norm{X^{(m)}-X^{(n)}}^{2} \right) 
& \leq & 4 \Exp \norm{ \sum_{k=n+1}^{m} \lambda_{k} \widetilde{M}_{T}(h_{k}) g_{k} }^{2} \\
& = & 4 \sum_{k=n+1}^{m} \abs{\lambda_{k}}^{2} \Exp \abs{ \widetilde{M}_{T}(h_{k})}^{2} \\
& \leq & 4  \norm{\widetilde{M}_{T}}_{\mathcal{L}(H, L^{2} \ProbSpace )}^{2} \sum_{k=n+1}^{m} \abs{\lambda_{k}}^{2}  \rightarrow 0, \mbox{ as } m,n \rightarrow 0. 
\end{eqnarray*} 
Then the sequence $(X^{(n)}:n \in \N)$ has a limit 
$$  X= \sum_{k=1}^{\infty} \lambda_{k} \widetilde{M}(h_{k}) g_{k} \in \mathcal{M}^{2}_{T}(G). $$
Furthermore for each $g \in G$, from the continuity of $\widetilde{M}$ we have $\Prob$-a.e
\begin{equation*}
\inner{X}{g}=  \sum_{k=1}^{\infty} \lambda_{k} \widetilde{M}(h_{k}) \inner{g_{k}}{g} = \widetilde{M} \left( \sum_{k=1}^{\infty} \lambda_{k} \inner{g_{k}}{g} h_{k} \right)=  \widetilde{M} (S^{*}g).    
\end{equation*}
Thus we conclude that $X$ is a version of $M \circ S^{*}$. 
\end{prf}

\section{Cylindrical Martingale-Valued Measures} \label{SectCylMarValMeasu}

The class of integrators for our stochastic integrals is the familly of cylindrical martingale-valued measures in the separable Hilbert space $H$. We recall its definition from \cite{FonsecaMora:2018-1}.

\begin{defi} \label{defiMartingValuedMeasure} 
Let $U$ be a topological space and consider a ring $\mathcal{R}\subseteq \mathcal{B}(U)$ that generates $\mathcal{B}(U)$.  A \emph{cylindrical martingale-valued measure} on $\R_{+} \times \mathcal{R}$ is a collection $M=(M(t,A): t \geq 0, A \in \mathcal{R})$ of cylindrical random variables in $H$ such that:
\begin{enumerate}
\item $\forall \, A \in \mathcal{R}$, $M(0,A)(h)= 0$ $\Prob$-a.e., $\forall h \in H$.
\item $\forall t \geq 0$, $M(t,\emptyset)(h)= 0$ $\Prob$-a.e. $\forall h \in H$ and if $A, B \in \mathcal{R}$ are disjoint then $\forall h \in H$, 
$$M(t,A \cup B)(h)= M(t,A)(h) + M(t,B)(h) \quad  \Prob \mbox{-a.e.}.$$
\item $\forall \, A \in \mathcal{R}$, $(M(t,A): t \geq 0)$ is a cylindrical mean-zero square integrable martingale and for each $t \geq 0$ the map $ M(t,A): H \rightarrow L^{0} \ProbSpace$ is continuous. 
\item For disjoint $A, B \in \mathcal{R}$, $\Exp \left( M(t,A)(h_{1}) M(s,B)(h_{2}) \right)=0$, for each $t,s \geq 0$, $h_{1}, h_{2} \in H$. 
\end{enumerate}
\end{defi}

We will further assume that the following properties are satisfied:
\begin{enumerate} \setcounter{enumi}{4}
\item For $0\leq s < t$, $M ( (s, t], A)(h) \defeq (M(t,A)- M(s,A))(h)$ is independent of $\mathcal{F}_{s}$, for all $A \in \mathcal{R}$, $h \in H$. 
\item For each $A \in \mathcal{R}$ and $0 \leq s < t$, 
\begin{equation} \label{covarianceFunctionalNuclearMartValuedMeasure}
\Exp \left( \abs{ M((s,t],A)(h)}^{2} \right) = \int_{s}^{t} \int_{A} q_{r,u}(h)^{2} \mu(du) \lambda (dr) , \quad \forall \, h \in H.
\end{equation} 
where 
\begin{enumerate}
	\item $\mu$ is a $\sigma$-finite measure on $(U, \mathcal{B}(U))$ satisfying $\mu(A)< \infty$, $\forall \, A \in \mathcal{R}$,
	\item $\lambda$ is a $\sigma$-finite measure on $(\R_{+},\mathcal{B}(\R_{+}))$, finite on bounded intervals,
	\item $\{q_{r,u}: r \in \R_{+}, \, u \in U \}$ is a family of  continuous Hilbertian semi-norms on $H$, such that for each $h_{1}$, $h_{2}$ in $H$, the map $(r,u) \mapsto q_{r,u}(h_{1},h_{2})$ is $\mathcal{B}(\R_{+}) \otimes \mathcal{B}(U)/ \mathcal{B}(\R_{+})$-measurable and  bounded on $[0,T] \times U$ for all $T>0$. Here, $q_{r,u}(\cdot,\cdot)$ denotes the positive, symmetric, bilinear form associated to the Hilbertian semi-norm $q_{r,u}$.  
\end{enumerate}
\end{enumerate}

\begin{exam} \label{examLevyMartValuedMeasure}
Let $L=(L_{t}: t \geq 0)$ be a $H$-valued c\`{a}dl\`{a}g L\'{e}vy process, i.e. $L$ has independent and stationary increments and has $\Prob$-a.e. c\`{a}dl\`{a}g paths. Assume furthermore that $L$ is $(\mathcal{F}_{t})$-adapted and that $L_{t}-L_{s}$ is independent of $\mathcal{F}_{s}$ for all $0 \leq s < t$. We can associate to $L$ a Poisson random measure $N$ on $\R^{+} \times (H \setminus \{ 0 \})$ given by 
$$ N(t,A)= \# \{ 0 \leq s  \leq t: \Delta L_{t} \defeq L_{t}-L_{t-} \in A \}. $$ 
Denote by $\mathcal{A}$ the collection of all $A \in \mathcal{B}(H \setminus \{ 0 \})$ such that $0 \notin \overline{A}$. It is easy to check that $\mathcal{A}$ is a ring. Moreover, it is well-known that for all $A \in \mathcal{A}$, $(N(t,A): t \geq 0)$ is a Poisson process with $\Exp \left( N(t,A) \right)= t \nu (A)$, where $\nu$ is a L\'{e}vy measure, i.e. $\nu$ is a Borel measure on $H$ with $\nu( \{ 0 \})=0$, and  $\displaystyle{ \int_{H} \norm{h}^{2} \wedge 1 \,  \nu (dh)< \infty}$ (see e.g. Chapter 2 in \cite{ApplebaumLPSC} for the finite dimensional case). 

Let $A \in \mathcal{A}$. For $f: A \rightarrow H$ measurable, we may define the Poisson integral as the finite random sum:
$$ \int_{A} f(h) \, N (t, dh) = \sum_{0 \leq s \leq t} f (\Delta L_{s}) \ind{A}{\Delta L_{s}}. $$ 
Moreover, let $\widetilde{N}(dt,dh)=N(dt,dh)-dt \nu(dh)$ be the compensated Poisson random measure corresponding to $N$. For $f \in L^{2}\left(A, \restr{\nu}{A}; H \right)$ we define the compensated Poisson integral:
$$ \int_{A} f(h) \, \widetilde{N} (t, dh) =\int_{A} f(h) \, N (t, dh)-t \int_{A} f(h) \nu(dh). $$
It is well-known (see e.g. Chapter 2 in \cite{ApplebaumLPSC}) that $\displaystyle{ \left( \int_{A} f(h) \widetilde{N}(t,dh): t \geq 0 \right)}$ is a mean-zero square integrable $H$-valued c\`{a}dl\`{a}g martingale, and 
$$ \Exp \left( \norm{ \int_{A} f(h) \widetilde{N}(t,dh) }^{2} \right)= t \int_{A} \norm{f(h)}^{2} \nu(dh). $$  
Recall that $L$ being a L\'{e}vy process in a separable Hilbert space, possesses a L\'{e}vy-It\^{o} decomposition (see \cite{AlbeverioRudiger:2005} for a proof, there in more general context of type 2 Banach spaces):
\begin{equation}\label{eqLevyItoDecomp}
L_{t}= t \xi +W_{t}+\int_{\norm{h} <1} \, h \,  \widetilde{N}(t,dh)+ \int_{\norm{h} \geq 1} \, h \, N(t,dh),  
\end{equation}
where $W=(W_{t}: t \geq 0)$ is a Wiener process in $H$ with covariance operator $Q$ (i.e.  $\Exp ( (h,W_{t})^{2})=t \inner{h}{Qh}$) which is positive and of trace class, $W$ is independent of the Poisson random measure $N$, $\left( \int_{\norm{h} \geq 1} \, h \, N(t,dh): t \geq 0 \right)$ is a Poisson integral,  and $\left( \int_{\norm{h} <1} \, h \,  \widetilde{N}(t,dh): t \geq 0 \right)$ is a $H$-valued c\`{a}dl\`{a}g mean-zero square integrable martingale such that 
$$ \int_{\norm{h} <1} \, h \,  \widetilde{N}(t,dh) = \lim_{n \rightarrow \infty} \int_{\frac{1}{n} \leq \norm{h} <1} \, h \,  \widetilde{N}(t,dh) \quad \forall t >0 \quad (\mbox{limit in }  L^{2}\ProbSpace ). $$ 

We can associate to the L\'{e}vy process $L$ a $H$-valued martingale-valued measure in the following way. Let $U \in \mathcal{B}(H)$ be such that $0 \in U$ and $\int_{U} \, \norm{h}^{2} \nu(du)< \infty$. Take $\mathcal{R}= \{  U \cap \Gamma: \Gamma \in \mathcal{A} \} \cup \{ \{0\}\}$, and let $M=(M(t,A): r \geq 0, A \in \mathcal{R})$ be given by
\begin{equation} \label{levyMartValuedMeasExam} 
M(t,A) = W_{t} \delta_{0}(A) + \int_{A \setminus \{0 \}} \,  h \, \widetilde{N}(t,dh), \quad \forall \, t \geq 0, \, A \in \mathcal{R}. 
\end{equation}
Then, $M$ is a $H$-valued martingale valued measure as defined for example in \cite{Applebaum:2006}. The induced cylindrical process $h \mapsto M(t,A)(h) \defeq \inner{M(t,A)}{h}$
 is a cylindrical martingale-valued measure, where for each $h \in H$, $0 \leq s \leq t$, $A \in \mathcal{R}$, 
\begin{equation} \label{secondMomentLevyNuclearMartValuedMeasu}
\Exp \left( \abs{M((s,t],A)(h)}^{2} \right)=(t-s) \left[ \inner{h}{Qh} + \int_{A \setminus \{0 \}} \inner{u}{h}^{2}  \, \nu(du) \right]. 
\end{equation}
In particular, with respect to the notation in Definition \ref{defiMartingValuedMeasure} we have: \begin{inparaenum}[(i)]
\item $\mu = \delta_{0}+ \restr{\nu}{U}$,  
\item $\lambda$ is the Lebesgue measure on $(\R_{+}, \mathcal{B}(\R_{+}))$, 
\item $\{q_{r,u}: r \in \R_{+}, u \in U \}$ is given by 
\begin{equation}\label{defiSemiNormsLevyMartValuedMeas}
q_{r,u}(h)^{2}= \begin{cases} \inner{h}{Qh}, & \mbox{if } u=0, \\  \inner{u}{h}^{2}, & \mbox{if } u \in U \setminus \{0\}. \end{cases}
\end{equation}
\end{inparaenum}
Following \cite{FonsecaMora:2018-1}, we call $M$ defined in \eqref{levyMartValuedMeasExam} a \emph{L\'{e}vy martingale-valued measure}. Indeed, $M$ also defines a nuclear decomposable martingale-valued measure as defined in \cite{Applebaum:2006}.  
\end{exam}

\begin{exam}\label{examCylLevySquareMarting}
Following \cite{ApplebaumRiedle:2010}, a cylindrical process $L=(L_{t}: t \geq 0)$ in $H$ is called a \emph{cylindrical L\'{e}vy process} if for all $n \in \N$ and $h_{1}, \dots, h_{n} \in H$, the $\R^{n}$-valued process $(L_{t}(h_{1}),\dots, L_{t}(h_{n}): t \geq 0)$ is a L\'{e}vy process. Assume further that $L$ has mean-zero and is square integrable, then $L$ is a cylindrical mean-zero square integrable martingale. Moreover, it follows from Corollary 3.12 in \cite{ApplebaumRiedle:2010}, that $L$ can be decomposed into $L_{t}=W_{t}+P_{t}$. Here, $W=(W_{t}:t \geq 0)$ is a cylindrical Wiener process, i.e. for all $n \in \N$ and $h_{1}, \dots, h_{n} \in H$, the $\R^{n}$-valued process $(W_{t}(h_{1}),\dots, W_{t}(h_{n}): t \geq 0)$ is a Wiener process. In particular, associated to $W$ is a quadratic functional $q:H \mapsto \R_{+}$ satisfying for $0 \leq s \leq t$ that 
$$ \Exp (\abs{ (W_{t}-W_{s})(h)}^{2})=(t-s) q(h), \quad \forall \, h \in H, \, t \geq 0.$$ 
The cylindrical process $P=(P_{t}: t \geq 0)$ is independent of $W$, and is a cylindrical mean-zero  square integrable martingale given by 
$$ P_{t}(h)= \int_{\R \setminus \{0\}} \, \beta \widetilde{N}_{h} (t, d\beta), \quad \forall \, h \in H, $$
where $ N_{h}$ denotes the Poisson random measure associated with the L\'{e}vy process $(P_{t}(h): t \geq 0)$ and with corresponding L\'{e}vy measure $\nu_{h}$. Indeed, Theorem 2.7 in \cite{ApplebaumRiedle:2010} shows that there exists a cylindrical measure $\nu$ on $\mathcal{Z}(H)$, called the cylindrical L\'{e}vy measure of $L$, such that for each $h \in H$ we have $\nu_{h}=\nu \circ \pi_{h}^{-1}$. In particular, for $0 \leq s \leq t$ we have
$$ \Exp (\abs{ (P_{t}-P_{s})(h)}^{2})=(t-s) \int_{H} \inner{y}{h}^{2} \nu(dy), \quad \forall \, h \in H, \, t \geq 0.$$  
We will assume that for each $t \geq 0$, the maps $W_{t}:H \rightarrow L^{0}\ProbSpace$ and $P_{t}:H \rightarrow L^{0}\ProbSpace$ are continuous (hence $L_{t}:H \rightarrow L^{0}\ProbSpace$ is continuous). Moreover, the continuity of $W_{1}$ implies the continuity of the quadratic form $q$. Hence, there exists $Q \in \mathcal{L}(H,H)$ satisfying $q(h)=\inner{h}{Qh}$ for each $h \in H$. Under the above assumptions, we can define a cylindrical martingale-valued measure by means of the prescription:
 \begin{equation} \label{eqCylLevySquare}
M(t,A)=(W_{t}+P_{t}) \delta_{0}(A), \quad \forall \, t \in \R_{+}, \, A \in \mathcal{B}(\{ 0 \}).  
\end{equation}
 Moreover, for each $0 \leq s \leq t$, we have:
\begin{equation} \label{eqCovariCylMartMeasLevySquar}
\Exp \left( \abs{ M((s,t],\{ 0 \})(h)}^{2} \right) = (t-s) \left[ \inner{h}{Qh}+ \int_{H} \inner{y}{h}^{2} \nu(dy) \right], \quad \forall \, h \in H.
\end{equation} 
With respect to the notation in Definition \ref{defiMartingValuedMeasure} we have: \begin{inparaenum}[(i)]
\item $U=\{0\}$, $\mathcal{R}=\mathcal{B}(\{ 0 \})$, 
\item $\mu = \delta_{0}$, $\lambda$ is the Lebesgue measure on $(\R_{+}, \mathcal{B}(\R_{+}))$, 
\item $\displaystyle{q_{r,0}(h)^{2}= \inner{h}{Qh}+ \int_{H} \inner{y}{h}^{2} \nu(dy)}$ $\forall \, r \geq 0$. 
\end{inparaenum}
\end{exam}
     
\begin{exam}\label{examLevyValRandMartMeasu}
Let $\mathcal{O} \subseteq \R^{d}$ be a Borel set and consider the $\delta$-ring $\mathcal{B}_{b}(\mathcal{O})=\{ A \in \mathcal{B}(\mathcal{O}): A \mbox{ is relatively compact}\}$. Following \cite{GriffithsRiedle:2019}, a \emph{L\'{e}vy-valued random martingale measure on $\mathcal{B}_{b}(\mathcal{O})$} is a family $(M(t,A): t \geq 0, A \in \mathcal{B}_{b}(\mathcal{O})) \subseteq L^{0} \ProbSpace$ such that 
\begin{enumerate}
\item For each $t \geq 0$, and for each collection of disjoint sets $A_{1}, A_{2}, \dots \in \mathcal{B}_{b}(\mathcal{O})$ the following hold:
\begin{enumerate}
\item The random variables $M(t,A_{1})$, $M(t,A_{2})$, $\dots$ are independent. 
\item If $\bigcup_{k \in \N} A_{k} \in \mathcal{B}_{b}(\mathcal{O})$ then $M \left(  \bigcup_{k \in \N} A_{k} \right) = \sum_{k \in \N} M(A_{k})$ $\Prob$-a.e.
\item The random variable $M(t,A)$ is infinitely divisible for each $A \in \mathcal{B}_{b}(\mathcal{O})$. 
\end{enumerate}
\item For every $A_{1}, A_{2}, \dots, A_{n} \in \mathcal{B}_{b}(\mathcal{O})$ and $n \in \N$, the process $((M(t,A_{1}), \dots, M(t,A_{n}): t \geq 0)$ is a L\'{e}vy process in $\R^{n}$. 
\end{enumerate} 
If we further assume in (2) that $((M(t,A_{1}), \dots, M(t,A_{n}): t \geq 0)$ is a zero-mean square integrable L\'{e}vy process, then $M$ defines a martingale-valued measure with values in $\R$ which satisfies all the conditions in Definition \ref{defiMartingValuedMeasure} for $U=\mathcal{O}$ and $\mathcal{R}=\mathcal{B}_{b}(\mathcal{O})$. In particular, as pointed out in Remark 2.4 in \cite{GriffithsRiedle:2019}, the classical definition of Gaussian space-time white noise fits naturally in the setting of the L\'{e}vy-valued random martingale measures, and hence also satisfies all the conditions in Definition \ref{defiMartingValuedMeasure}.  

As shown in \cite{GriffithsRiedle:2019}, the above concept can be related with that of cylindrical L\'{e}vy processes. Indeed, for any given locally finite Borel measure $\zeta$ on $\mathcal{O}$, let $H=L^{2}(\mathcal{O},\mathcal{B}(\mathcal{O}),\zeta)$ and consider a cylindrical process $L=(L_{t}: t \geq 0)$ in $H$. Following \cite{GriffithsRiedle:2019}, we say that $L$  is \emph{independently scattered} if for every collection of disjoint sets $A_{1}, A_{2}, \dots, A_{n} \in \mathcal{B}_{b}(\mathcal{O})$  and $n \in \N$, the random variables $L(t) \mathbbm{1}_{A_{1}}, \dots, L(t) \mathbbm{1}_{A_{n}}$ are independent for each $t \geq 0$. In such case, Theorem 4.5 in  \cite{GriffithsRiedle:2019} shows that the prescription 
$$M(t,A) \defeq L(t) \mathbbm{1}_{A}, \quad \forall \, t \geq 0, \, A \in \mathcal{B}_{b}(\mathcal{O}), $$ 
defines a L\'{e}vy-valued random martingale measure on $\mathcal{B}_{b}(\mathcal{O})$. If $L$ has zero-mean and second moments, we again obtain an example that satisfies the conditions in Definition \ref{defiMartingValuedMeasure} (but which is different from  Example \ref{examCylLevySquareMarting}). 

As an example we have the impulsive cylindrical process in $L^{2}(\mathcal{O},\mathcal{B}(\mathcal{O}),\zeta)$ defined by Peszat and Zabczyk in Section 7.2 in \cite{PeszatZabczyk} (see also Example 3.6 in \cite{ApplebaumRiedle:2010}), and given by 
$$ L(t)f \defeq \int_{0}^{t} \int_{\mathcal{O}}\int_{\R} f(x) y \widetilde{N}(ds,dx,dy), $$
where $N$ is a Poisson random measure on $\R_{+} \times \mathcal{O} \times \R$ with intensity $\mbox{leb} \otimes \zeta \otimes \mu$ for a L\'{e}vy measure $\mu$ on $\mathcal{B}(\R)$ satisfying $\int_{\R} \beta^{2} \mu(d\beta)<\infty$. This cylindrical process is independently scattered as shown in Example 4.8 in  \cite{GriffithsRiedle:2019}, has mean-zero and second moment.   
\end{exam}

We finalize this section with the following technical result on the family of Hilbertian seminorms $\{q_{r,u}: r \in \R_{+}, \, u \in U \}$ in Definition \ref{defiMartingValuedMeasure}(6). Its conclusion will be of importance for our theory of stochastic integration. 

\begin{lemm} \label{lemmInecuaForSeminorms} There exists  $C>0$ such that $q_{r,u}(\cdot) \leq C \norm{\cdot}$ for each $(r,u) \in [0,T] \times U$.
\end{lemm}
\begin{prf}
Following the proof of Theorem 4.2 in \cite{BojdeckiJakubowski:1989}, for each $(r,u) \in [0,T] \times U$ let $V_{r,u}=\{ h \in H: q_{r,u}(h) \leq 1 \}$ and $U=\bigcap_{(r,u)} V_{r,u}$. Since each $V_{r,u}$ is the closed unit ball of $q_{r,u}$, the set $U$ is closed, convex and balanced. Moreover, from Definition \ref{defiMartingValuedMeasure}(6)(c), we have that for each $h \in H$, $E(h) \defeq \sup\{ q_{r,u}(h): (r,u) \in [0,T] \times U \}< \infty$. Therefore, $h \in E(h) U$, showing that $U$ is absorbing, hence a barrel in $H$. But as $H$ is Hilbert (hence a barrelled space), we have that $U$ is a neighborhood of zero. Thus there exists some $C>0$ such that $\{ h \in H: \norm{h} \leq 1/C \} \subseteq U \subseteq V_{r,u}$. Therefore $q_{r,u}(\cdot) \leq C \norm{\cdot}$ for each $(r,u) \in [0,T] \times U$.
\end{prf}


\section{The Stochastic Integral}\label{sectStochaIntegral}

From now on $H$ and $G$ will denote two separable Hilbert spaces and $M$ is a  cylindrical martingale-valued measure in $H$ as in Definition \ref{defiMartingValuedMeasure}. Recall that $H_{q_{r,u}}$ is the Hilbert space determined on $H$ by the Hilbertian seminorm $q_{r,u}$ and $i_{q_{r,u}}: H \rightarrow H_{q_{r,u}}$ is the canonical inclusion (Sect. \ref{sectionPrelim}). 

\subsection{Construction of the Stochastic Integral} \label{sectConstrStochInte}

The main objective of this section is to define the stochastic integral with respect to $M$ for the following class of operator-valued maps: 

\begin{defi} \label{defiIntegrandsWithSquareMoments}
Let $\Lambda^{2}(M,T; H,G)$ denote the collection  of families $\Phi=(\Phi(r,\omega,u): r \in [0,T], \omega \in \Omega, u \in U )$ of operator-valued maps satisfying the following: 
\begin{enumerate}
\item $\Phi(r,\omega,u) \in \mathcal{L}_{2}(H_{q_{r,u}}, G)$, for all $r \in [0, T]$, $\omega \in \Omega$, $u \in U$, 
\item $\Phi$ is \emph{$q_{r,u}$-predictable}, i.e. for each $h \in H$, $g \in G$, the mapping $[0,T] \times \Omega \times U \rightarrow \R_{+}$ given by $(r,\omega,u) \mapsto q_{r,u}(\Phi(r,\omega,u)^*g, i_{q_{r,u}} h)$ is $\mathcal{P}_{T} \otimes \mathcal{B}(U)/\mathcal{B}(\R_{+})$-measurable.
\item 
\begin{equation} \label{eqDefiNormSquaIntegra}
\norm{\Phi}_{\Lambda^{2}(M,T)}^{2} \defeq \Exp \int_{0}^{T} \int_{U}  \norm{ \Phi(r,u)}_{\mathcal{L}_{2}(H_{q_{r,u}}, G)}^{2} \mu(du) \lambda(dr) < \infty.
\end{equation} 
\end{enumerate}
We will write $ \Lambda^{2}(M,T;H,G)$ by $\Lambda^{2}(M,T)$ for simplicity when the underlying spaces are clear.  
\end{defi}

\begin{rema} \label{remaNormOfIntegrandsIsPredictable}
Since the Hilbert space $H$ is separable, it follows from Definition \ref{defiIntegrandsWithSquareMoments}(2), and from Proposition 3.8 in \cite{FonsecaMora:2018-1}, that for each $g \in G$, the map  $(r,\omega,u) \mapsto q_{r,u}(\Phi(r,\omega,u)^*g)^{2}$ is $\mathcal{P}_{T} \otimes \mathcal{B}(U)/\mathcal{B}(\R_{+})$-measurable. Therefore, because $G$ is a separable Hilbert space, the map $(r,\omega,u) \mapsto \norm{ \Phi(r,\omega, u)}_{\mathcal{L}_{2}(H_{q_{r,u}}, G)}^{2}$ is $\mathcal{P}_{T} \otimes \mathcal{B}(U)/\mathcal{B}(\R_{+})$-measurable and the integrand in \eqref{eqDefiNormSquaIntegra} is well-defined.
\end{rema}

The proof of the following result can be carried out from similar arguments to those in the proof of Proposition 2.4 in \cite{BojdeckiJakubowski:1990}. 

\begin{prop} \label{propSpaceSquareIntegHilbertSpace} $ \Lambda^{2}(M,T)$ is a Hilbert space when equipped with the inner product $\inner{\cdot}{\cdot}_{\Lambda^{2}(M,T)}$ corresponding to the Hilbertian norm $\norm{\cdot }_{\Lambda^{2}(M,T)}$.
\end{prop}

Our next objective is to define the stochastic integral for the integrands that belongs to the space $ \Lambda^{2}(M,T)$. To do this, we will follows the standard approach to define the integral first for a subclass of simpler processes that we define as follows. 

\begin{defi} \label{defiClassSimpleInteg}
Let $\mathcal{S}(M,T)$ be the collection of all the families $\Phi=\{ \Phi(r,\omega,u): r \in [0,T],  \omega \in \Omega, u \in U\}$ of operator-valued maps of the form:
\begin{equation} \label{eqSimpleIntegrands}
\Phi(r,\omega,u)=  \sum_{i=1}^{n}\sum_{j=1}^{m} \ind{]s_{i}, t_{i}]}{r} \ind{F_{i}}{\omega} \ind{A_{j}}{u} S_{i,j} \circ i^{*}_{q_{r,u}},  
\end{equation}
for all $r \in [0,T]$, $\omega \in \Omega$, $u \in U$, where $m$, $n \in \N$, and for $i=1, \dots, n$, $j=1, \dots, m$, $0\leq s_{i}<t_{i} \leq T$, $F_{i} \in \mathcal{F}_{s_{i}}$, $A_{j} \in \mathcal{R}$ and $S_{i,j} \in \mathcal{L}_{2}(H,G)$.
\end{defi}

\begin{rema}
Let $\Phi \in \mathcal{S}(M,T)$ be of the form \eqref{eqSimpleIntegrands}. Observe that (taking a smaller partition if needed) we can always assume that the following is satisfied: 
\begin{equation} \label{separationPartitionIntervalsSimpleForm}
\mbox{for } k \neq j, \quad  ]s_{k},t_{k}]  \, \cap  \,  ]s_{j},t_{j}] \neq \emptyset \quad \Rightarrow  \quad ]s_{k},t_{k}]= \, ]s_{j},t_{j}] \mbox{ and } F_{k} \cap F_{j}= \emptyset.
\end{equation}  
\end{rema}

It is easy to check that $\mathcal{S}(M,T) \subseteq  \Lambda^{2}(M,T)$. In particular, observe that since by Lemma \ref{lemmInecuaForSeminorms} we have $\sup_{(r,u)} \norm{i_{q_{r,u}}} \leq C$ for some $C>0$, then if $\Phi \in \mathcal{S}(M,T)$ is of the form \eqref{eqSimpleIntegrands},  we have 
$$ \norm{\Phi}_{\Lambda^{2}(M,T)}^{2} 
\leq C^{2} \sum_{i=1}^{n}\sum_{j=1}^{m} (t_{i}- s_{i}) \Prob(F_{i})\mu(A_{j}) \norm{S_{i,j}}^{2}_{\mathcal{L}_{2}(H,G)}  < \infty. $$
Moreover we have the following: 

\begin{prop} \label{propSimpleDenseInSquareIntegra}
$S(M,T)$ is dense in $\Lambda^{2}(M,T)$. 
\end{prop}
\begin{prf}
Let $C(M,T)$ be the collection of all families of operator-valued maps $\Psi=\{ \Psi(r,\omega,u): r \in [0,T],  \omega \in \Omega, u \in U\}$ taking the simple form
\begin{equation} \label{simpleFamiliesProofDenseSupspaceIntegrands}
\Psi(r,\omega,u)= \ind{]s,t]}{r} \ind{F}{\omega} \ind{A}{u} S \circ i^{*}_{q_{r,u}}, \quad \forall \, r \in [0,T], \, \omega \in \Omega, \, u \in U,
\end{equation}  
where $0 \leq s < t \leq T$, $F \in \mathcal{F}_{s}$, $A \in \mathcal{R}$ and $S \in \mathcal{L}_{2}(H,G)$. 

Since $C(M,T)$ spans $S(M,T)$, the orthogonal complements of these spaces satisfy $\overline{S(M,T)}^{\perp}=S(M,T)^{\perp}=C(M,T)^{\perp}$. Hence, to show that $S(M,T)$ is dense in $\Lambda^{2}(M,T)$ it is sufficient to show that $C(M,T)^{\perp}=\{0\}$

In effect, let $\Phi \in C(M,T)^{\perp}$. If $\Psi \in C(M,T)$ is of the form \eqref{simpleFamiliesProofDenseSupspaceIntegrands}, then we have that 
\begin{equation} \label{innerProductProofDenseSubpaceIntegrands}
\inner{\Phi}{\Psi}_{\Lambda^{2}(M,T)}= \int_{F} \int_{s}^{t} \int_{A} \inner{ \Phi(r,u)}{ S \circ i^{*}_{q_{r,u}}}_{\mathcal{L}_{2}(H_{q_{r,u}}, G)}  \mu(du) \lambda(dr) \Prob(d \omega)=0.   
\end{equation}
Let $(h_{n}: n \in \N)$ and $(g_{k}: k \in \N)$ be complete orthonormal systems in $H$ and $G$ respectively. Then, the collection $(h_{n} \otimes g_{k}: n, k \in \N)$ is a complete orthonormal system in $\mathcal{L}_{2}(H,G)$.
Because $\mathcal{P}_{T} \otimes \mathcal{B}(U)$ is generated by the family of all subsets of $[0,T] \times \Omega \times U$ of the form $G=]s,t] \times F \times A$, where $0 \leq s < t \leq T$, $F \in \mathcal{F}_{s}$, $A \in \mathcal{R}$;  then \eqref{innerProductProofDenseSubpaceIntegrands} and Fubini's theorem imply that $\lambda \otimes \Prob \otimes \mu$-a.e. 
$$q_{r,u}(\Phi(r,u)^{*}g_{k}, i_{q_{r,u}} h_{n})= \inner{ \Phi(r,u)}{ (h_{n} \otimes g_{k})\circ i^{*}_{q_{r,u}}  }_{\mathcal{L}_{2}(H_{q_{r,u}}, G)}=0, \quad \forall n, k  \in \N.$$ 
Since the orthonormal systems $(h_{n}: n \in \N)$ and $(g_{k}: k \in \N)$ are complete and each $q_{r,u}$ is continuous, the above implies that 
$\lambda \otimes \Prob \otimes \mu$-a.e. 
\begin{equation}\label{eqWeakImageCeroAdjointSeminorms}
q_{r,u}(\Phi(r,u)^{*}g, i_{q_{r,u}} h)=0, \quad \forall g \in G, h \in H.
\end{equation}
But because $i_{q_{r,u}}(H)$ is dense in $H_{q_{r,u}}$,  \eqref{eqWeakImageCeroAdjointSeminorms} implies that $\lambda \otimes \Prob \otimes \mu$-a.e. we have $\Phi(r,u)^{*}g=0$ $\forall g \in G$, i.e. $\mbox{Kernel}(\Phi(r,u)^{*})=G$ and consequently $\mbox{Range}(\Phi(r,u))=\{0\}$. From the arguments given at the beginning of the proof, we conclude that  $S(M,T)$ is dense in $\Lambda^{2}(M,T)$.   
\end{prf}

Our next objective is to define the weak stochastic integral for the elements of $\mathcal{S}(M,T)$. To do this, we will need to set some notation. 

\begin{defi}
Given $0 \leq s < t \leq T$, $A \in \mathcal{R}$ and $S \in \mathcal{L}_{2}(H,G)$, we denote by $(Y_{r}(s,t,A,S): r \in [0,T])$ the $G$-valued c\`{a}dl\`{a}g square integrable martingale obtained via radonification of the cylindrical square integrable martingale $h \mapsto \left( M((s \wedge r, t \wedge r],A)(h): t \in [0,T] \right)$ by the Hilbert-Schmidt operator $S$ (see Theorem \ref{theoRadonCylMartHilbSchmi}), i.e. $(Y_{r}(s,t,A,S): r \in [0,T])$ satisfies that $\Prob$-a.e. 
\begin{equation}\label{eqDefRadonProcessY}
\inner{Y_{r}(s,t,A,S)}{g}= M((s \wedge r, t \wedge r],A)(S^{*}g), \quad \forall r \in [0,T], \, g \in G. 
\end{equation}
\end{defi}
  
We are ready to define our stochastic integral for simple integrands via radonification:

\begin{defi}
Let $\Phi \in \mathcal{S}(M,T;H,G)$ be of the form \eqref{eqSimpleIntegrands} and satisfying \eqref{separationPartitionIntervalsSimpleForm}. We define the \emph{stochastic integral of $\Phi$} as the $G$-valued process $(I_{t}(\Phi): t \in [0,T])$ defined by  
\begin{equation} \label{eqDefiStochasticIntegralSimpleIntegrands}
I_{t}(\Phi)= \sum_{i=1}^{n} \sum_{j=1}^{m} \mathbbm{1}_{F_{i}} Y_{t}(s_{i},t_{i},A_{j},S_{i,j}).
\end{equation}
We will write $ \mathcal{S}(M,T;H,G)$ by $\mathcal{S}(M,T)$ for simplicity when the underlying spaces are clear.  
\end{defi}

It should be clear to the reader that from the finite additivity of $M$ on $\mathcal{R}$, and from \eqref{eqDefRadonProcessY}, the process $(I_{t}(\Phi): t \in [0,T])$ is uniquely defined (up to indistinguishable versions) independently of the representation of $\Phi \in \mathcal{S}(M,T)$ (i.e. of the expression of $\Phi$ as in \eqref{eqSimpleIntegrands}). 

\begin{rema} \label{remaCyliAndClassicIntegCoincide}
Suppose that $M$ is a $H$-valued martingale-valued measure (e.g. the L\'{e}vy martingale-valued measure in Example \ref{examLevyMartValuedMeasure}). In that case, \eqref{eqDefRadonProcessY} takes the form 
$$   \inner{Y_{r}(s,t,A,S)}{g}= \inner{S \circ M((s \wedge r, t \wedge r],A)}{g}, \quad \forall r \in [0,T], \, g \in G. $$
Hence, if $\Phi \in \mathcal{S}(M,T;H,G)$ is of the form \eqref{eqSimpleIntegrands} and satisfying \eqref{separationPartitionIntervalsSimpleForm}, then \eqref{eqDefiStochasticIntegralSimpleIntegrands} yields: 
$$I_{t}(\Phi)= \sum_{i=1}^{n} \sum_{j=1}^{m} \mathbbm{1}_{F_{i}} S_{i,j} \circ  M((s_{i} \wedge t,t_{i}\wedge t ],A_{j}).$$
Therefore, our construction coincides with the classical definition of stochastic integral in Hilbert spaces as for example in \cite{Applebaum:2006, DaPratoZabczyk, PeszatZabczyk}. 
\end{rema}

The basic properties of the stochastic integral are given in the next result. Recall from Sect. \ref{sectionPrelim} that $\mathcal{M}_{T}^{2}(G)$ denotes the space of $G$-valued mean-zero square integrable c\`{a}dl\`{a}g martingales.

\begin{theo} \label{theoItoIsometrySimpleIntegrands}
For every $\Phi \in \mathcal{S}(M,T)$, $(I_{t}(\Phi): t \in [0,T]) \in \mathcal{M}_{T}^{2}(G)$ and for each $t \in [0,T]$,        
\begin{equation}\label{eqItoIsometrySimpleIntegrands}
\Exp (I_{t}(\Phi))=0, \quad 
\Exp \left( \norm{I_{t}(\Phi)}^{2} \right) = \Exp \int_{0}^{t} \int_{U}  \norm{ \Phi(r,u)}_{\mathcal{L}_{2}(H_{q_{r,u}}, G)}^{2} \mu(du) \lambda(dr). 
\end{equation}
In particular, the map $I: \mathcal{S}(M,T) \rightarrow \mathcal{M}^{2}_{T}(G)$, $\Phi \mapsto (I_{t}(\Phi): t \in [0,T])$ is linear continuous.  
\end{theo}
\begin{prf}
Let $\Phi \in \mathcal{S}(M,T)$ be of the form \eqref{eqSimpleIntegrands} and satisfying \eqref{separationPartitionIntervalsSimpleForm}. The fact that  $(I_{t}(\Phi): t \in [0,T]) \in \mathcal{M}_{T}^{2}(G)$ is a consequence of the fact that the process $(Y_{r}(s,t,A,S): r \in [0,T])$ satisfies these properties and from \eqref{eqDefiStochasticIntegralSimpleIntegrands}. 

To prove \eqref{eqItoIsometrySimpleIntegrands}, let $(g_{d}: d \in \N)$ be a complete orthonormal system in $G$. From \eqref{eqDefRadonProcessY} and \eqref{eqDefiStochasticIntegralSimpleIntegrands} we have 
\begin{flalign*} 
& \Exp \left( \norm{I_{t}(\Phi)}^{2} \right) \\
& =  \Exp \sum_{d=1}^{\infty} \abs{ \inner{I_{t}(\Phi)}{g_{d}}}^{2} \\ 
& =  \Exp \sum_{d=1}^{\infty} \abs{ \sum_{i=1}^{n} \sum_{j=1}^{m} \mathbbm{1}_{F_{i}} \inner{Y_{t}(s_{i},t_{i},A_{j},S_{i,j})}{g_{d}} }^{2} \\
& =  \sum_{d=1}^{\infty} \Exp  \abs{ \sum_{i=1}^{n} \sum_{j=1}^{m} \mathbbm{1}_{F_{i}}  M((s_{i} \wedge t, t_{i} \wedge t],A_{j})(S_{i,j}^{*} g_{d}) }^{2} \\
& =  \sum_{d=1}^{\infty}  \sum_{i,k=1}^{n} \sum_{j,l=1}^{m} \Exp  \left( \mathbbm{1}_{F_{i}}  M((s_{i} \wedge t, t_{i} \wedge t ],A_{j})(S_{i,j}^{*} g_{d}) \cdot \mathbbm{1}_{F_{k}} M((s_{l} \wedge t , t_{l} \wedge t ],A_{k})(S_{k,l}^{*} g_{d})  \right)
\end{flalign*}
Now, from the properties of orthogonality of $M$ on the ring $\mathcal{R}$ (see Definition \ref{defiMartingValuedMeasure}(4)), we have that 
$$ \Exp \left( M((s_{i} \wedge t, t_{i} \wedge t ],A_{j})(S_{i,j}^{*} g_{d}) \cdot M((s_{l} \wedge t , t_{l} \wedge t ],A_{k})(S_{k,l}^{*} g_{d})  \right) = 0, $$
for each $d \in \N$, $i,k=1, \dots, n$, $j,l=1, \dots, m$, $i \neq k$, $j \neq l$. Hence,   
\begin{eqnarray*} 
\Exp \left( \norm{I_{t}(\Phi)}^{2} \right) 
& = & \sum_{d=1}^{\infty}  \sum_{i=1}^{n} \sum_{j=1}^{m} \Prob (F_{i}) \Exp \left( \abs{ M((s_{i} \wedge t, t_{i} \wedge t],A_{j})(S_{i,j}^{*}g_{d}) }^{2} \right)  \\
& = & \sum_{i=1}^{n} \sum_{j=1}^{m} \Prob (F_{i}) \sum_{d=1}^{\infty} \int_{s_{i} \wedge t}^{t_{i} \wedge t} \int_{A_{j}} q_{r,u}(i_{q_{r,u}} \circ S_{i,j}^{*} g_{d} )^{2} \mu (du) \lambda (dr)  \\
& = & \sum_{i=1}^{n} \sum_{j=1}^{m} \Prob (F_{i}) \int_{s_{i} \wedge t}^{t_{i} \wedge t} \int_{A_{j}} \norm{i_{q_{r,u}} \circ S_{i,j}^{*}}^{2}_{\mathcal{L}_{2}(G, H_{q_{r,u}})} \mu (du) \lambda (dr)  \\
& = & \Exp \int_{0}^{t} \int_{U}  \norm{ \Phi(r,u) }_{\mathcal{L}_{2}(H_{q_{r,u}}, G)}^{2} \mu(du) \lambda(dr).  
\end{eqnarray*}
It is easy to see from \eqref{eqDefiStochasticIntegralSimpleIntegrands} that the map $I: \mathcal{S}(M,T) \rightarrow \mathcal{M}^{2}(G)$ is linear. To prove the continuity, observe that from Doob's inequality,  \eqref{eqDefiNormSquaIntegra} and  \eqref{eqItoIsometrySimpleIntegrands}, for each $\Phi \in \mathcal{S}(M,T)$ we have
\begin{equation} \label{eqIntegralMapIsContinuousInequ}
\norm{I(\Phi)}^{2}_{\mathcal{M}^{2}_{T}(G)} = \Exp \left( \sup_{t\in [0,T]} \norm{I_{t}(\Phi)}^{2} \right) \leq 4 T \, \Exp \left( \norm{I_{T}(\Phi)}^{2} \right)= 4T \, \norm{\Phi}^{2}_{\Lambda^{2}(M,T)},  
\end{equation}   
The above inequality implies the continuity of the operator $I: \mathcal{S}(M,T) \rightarrow \mathcal{M}^{2}(G)$. 
\end{prf}

We can now extend the stochastic integral to integrands in $\Lambda^{2}(M,T)$. 
 
\begin{theo} \label{theoExisIntegralSquareInteg}
The map $I: \mathcal{S}(M,T) \rightarrow \mathcal{M}^{2}_{T}(G)$, $\Phi \mapsto (I_{t}(\Phi): t \in [0,T])$ have a continuous and linear extension $I: \Lambda^{2}(M,T) \rightarrow \mathcal{M}^{2}_{T}(G)$, such that for each $\Phi \in \Lambda^{2}(M,T)$ and each $t \in [0,T]$,    
\begin{equation}\label{eqItoIsometrySquareIntegrands}
\Exp (I_{t}(\Phi))=0, \quad 
\Exp \left( \norm{I_{t}(\Phi)}^{2} \right) = \Exp \int_{0}^{t} \int_{U}  \norm{ \Phi(r,u)}_{\mathcal{L}_{2}(H_{q_{r,u}}, G)}^{2} \mu(du) \lambda(dr). 
\end{equation}
\end{theo}	
\begin{prf} The existence of the continuous and linear extension $I: \Lambda^{2}(M,T) \rightarrow \mathcal{M}^{2}_{T}(G)$ is a direct consequence of Proposition \ref{propSimpleDenseInSquareIntegra} and Theorem \ref{theoItoIsometrySimpleIntegrands}. Now, for any given $t \in [0,T]$, it follows from \eqref{eqItoIsometrySimpleIntegrands} that the map $I_{t}: \mathcal{S}(M,T) \rightarrow L^{2}(\Omega, \mathcal{F}, \Prob; G)$, $\Phi \mapsto I_{t}(\Phi)$, is a linear isometry. Therefore, from Proposition \ref{propSimpleDenseInSquareIntegra} it has a linear isometric extension $I_{t}: \Lambda^{2}(M,T) \rightarrow L^{2}(\Omega, \mathcal{F}, \Prob; G)$. Then,  \eqref{eqItoIsometrySquareIntegrands} is a direct consequence of \eqref{eqItoIsometrySimpleIntegrands}.   
\end{prf}


\begin{defi}
We call the map $I$ defined in Theorem \ref{theoExisIntegralSquareInteg} the \emph{stochastic integral} mapping and for each $\Phi \in \Lambda^{2}(M,T)$, we call $I(\Phi)$ the \emph{stochastic integral} of $\Phi$. We will also denote the process $I(\Phi)$ by $\left\{ \int^{t}_{0} \int_{U} \Phi (r,u) M (dr, du): t \in [0,T] \right\}$.  
\end{defi}

\begin{prop}\label{propIntegContIfCylinMartIsCont}
If for each $A \in \mathcal{R}$ and $\phi \in \Phi$, the real-valued process $(M(t,A)(\phi): t \geq 0)$ is continuous, then for each $\Phi \in \Lambda^{2}(M,T)$ the stochastic integral $I(\Phi)$ is a continuous process.  
\end{prop}
\begin{prf}
For a simple process $\Phi \in \mathcal{S}(M,T)$, the result follows directly from \eqref{eqDefiStochasticIntegralSimpleIntegrands} and the fact that the radonified version $Y$ defined by \eqref{eqDefiStochasticIntegralSimpleIntegrands} is a continuous process if each $(M(t,A)(\phi): t \geq 0)$ is continuous. Since the sequential limit of $G$-valued continuous processes is continuous, the result now extends by the denseness of $\mathcal{S}(M,T)$ in $\Lambda^{2}(M,T)$ and the continuity of the stochastic integral mapping to every $\Phi \in \Lambda^{2}(M,T)$. 
\end{prf}

\subsection{Examples and Comparison with Other Theories of Integration} \label{sectExamStochInte}

This section is devoted to explore examples of our theory of stochastic integration with a special emphasis on the L\'{e}vy case. We will explore alternative characterizations for our integrands and compare our theory of integration with other theories available on the literature. To do this, it will be of great importance the following class of integrands. 

\begin{defi}
Let $\mathcal{H}^{2}(M,T)$ denote the collection of all mappings $\Psi: [0,T] \times \Omega \times U \rightarrow  \mathcal{L}_{2}(H, G)$ that are $\mathcal{P}_{T} \otimes \mathcal{B}(U)/\mathcal{B}(\mathcal{L}_{2}(H, G))$-measurable  and for which 
\begin{equation} \label{eqNormSpaceH2MT}
\norm{\Psi}_{\mathcal{H}^{2}(M,T)}^{2} \defeq \Exp \int_{0}^{T} \int_{U}  \norm{ \Psi(r,u)}_{\mathcal{L}_{2}(H, G)}^{2} \mu(du) \lambda(dr) < \infty.
\end{equation} 
\end{defi}

It is easy to check that $(\mathcal{H}^{2}(M,T), \inner{\cdot}{\cdot}_{\mathcal{H}^{2}(M,T)})$ is a Hilbert space. In the next result we will show that the elements in $\mathcal{H}^{2}(M,T)$ can be regarded as elements in $ \Lambda^{2}(M,T)$. 

\begin{prop} \label{propContEmbeddSquaIntegra} For every 
$\Psi \in \mathcal{H}^{2}(M,T)$, the family $\Psi \circ i^{*}_{q} \defeq (\Psi(r,\omega,u) \circ i^{*}_{q_{r,u}}: r \in [0,T], \omega\in \Omega, u \in U )$ is an element of $\Lambda^{2}(M,T)$.  Moreover, there exists  $C>0$ such that $\norm{\Psi \circ i^{*}_{q}}_{\Lambda^{2}(M,T)} \leq C \norm{\Psi}_{\mathcal{H}^{2}(M,T)}$ for every $\Psi \in \mathcal{H}^{2}(M,T)$.   
\end{prop}
\begin{prf}
Let $\Psi \in \mathcal{H}^{2}(M,T)$. It is clear that $\Psi(r,\omega,u) \circ i^{*}_{q_{r,u}} \in \mathcal{L}_{2}(H_{q_{r,u}},G)$ for each $r \in [0,T], \omega\in \Omega, u \in U$. Furthermore observe that from Definition \ref{defiMartingValuedMeasure}(6)(c) the mapping $(r,u) \mapsto q_{r,u}(i_{q_{r,u}}h_{1}, i_{q_{r,u}} h_{2})$ is $\mathcal{B}(\R_{+}) \otimes \mathcal{B}(U)/ \mathcal{B}(\R_{+})$-measurable for each $h_{1}, h_{2} \in H$, and since $\Psi$ is $\mathcal{P}_{T} \otimes \mathcal{B}(U)/\mathcal{B}(\mathcal{L}_{2}(H, G))$-measurable the map  
$$(r,\omega,u) \mapsto q_{r,u}( (\Psi(r,\omega,u) \circ i^{*}_{q_{r,u}})^*g, i_{q_{r,u}} h)= q_{r,u}(i_{q_{r,u}} \circ \Phi(r,\omega,u)^*g, i_{q_{r,u}} h), $$ 
is $\mathcal{P}_{T} \otimes \mathcal{B}(U)/\mathcal{B}(\R_{+})$-measurable for each $h \in H$, $g \in G$, i.e. $\Psi \circ i^{*}_{q}$ is $q_{r,u}$-predictable. 

Finally, observe that from Lemma \ref{lemmInecuaForSeminorms} there exists $C>0$ such that 
\begin{eqnarray*}
\norm{ \Psi(r,\omega,u) \circ i^{*}_{q_{r,u}} }_{\mathcal{L}_{2}(H_{q_{r,u}}, G)} 
&  = & \norm{ i_{q_{r,u}} \circ \Psi(r,u)^{*}}_{\mathcal{L}_{2}( G, H_{q_{r,u}})} \\
 & \leq & \sup_{(r,u)} \norm{ i_{q_{r,u}}}_{\mathcal{L}(H,H_{q_{r,u}})} \norm{ \Psi(r,u)^{*}}_{\mathcal{L}_{2}(G,H)}\\ 
& \leq & C \norm{ \Psi(r,u)^{*}}_{\mathcal{L}_{2}(G,H)}   
= C \norm{ \Psi(r,u)}_{\mathcal{L}_{2}(H, G)},   
\end{eqnarray*}
for each $(r,u) \in [0,T] \times U$. Therefore, the above inequality,  \eqref{eqDefiNormSquaIntegra} and \eqref{eqNormSpaceH2MT} show that we have 
$\norm{\Psi \circ i^{*}_{q}}_{\Lambda^{2}(M,T)} \leq C \norm{\Psi}_{\mathcal{H}^{2}(M,T)}$ and hence $\Psi \circ i^{*}_{q} \in \Lambda^{2}(M,T)$. 
\end{prf}

\begin{nota}
For now on, to every $\Psi \in \mathcal{H}^{2}(M,T)$ we will denote  $\norm{\Psi \circ i^{*}_{q} }_{\Lambda^{2}(M,T)}$ by $ \norm{\Psi}_{\Lambda^{2}(M,T)}$ and 
$\displaystyle{\int_{0}^{t}\int_{U} \Psi(r,u) \circ i^{*}_{q_{r,u}}  M(dr,du)}$ by $\displaystyle{\int_{0}^{t}\int_{U} \Psi(r,u) M(dr,du)}$. 
\end{nota}

The following result provides an alternative expression for the norm $\norm{\cdot }_{\Lambda^{2}(M,T)}$  for integrands in $\mathcal{H}^{2}(M,T)$. This formula will be of great importance in our examples. 

\begin{prop} \label{propCovarOperIntegrNorm}
Suppose that for each $r \geq 0$, $u \in U$, there exists a non-negative, symmetric, linear and continuous operator $Q_{r,u}:H \rightarrow H$, such that 
\begin{equation} \label{eqInnerSemiNormsCovaOper}
q_{r,u}(h_{1},h_{2})=\inner{h_{1}}{Q_{r,u}h_{2}}, \quad \forall \, h_{1},h_{2} \in H. 
\end{equation}
Then for every $\Psi \in \mathcal{H}^{2}(M,T)$ we have
\begin{equation} \label{eqEquivaIntegNormCovarOper}
 \norm{\Psi }_{\Lambda^{2}(M,T)}^{2}  =  \Exp \int_{0}^{T} \int_{U} \norm{\Psi(r,u) \circ Q^{1/2}_{q_{r,u}}}^{2}_{\mathcal{L}_{2}(H,G)} \mu(du) \lambda(dr). 
\end{equation}
\end{prop}
\begin{prf}
Let $(g_{k}: k \in \N)$ be a complete orthonormal system in $G$. Let $r \geq 0$, $\omega \in \Omega$, $u \in U$. Recall that for each $h \in H$, we have $q_{r,u}(i_{q_{r,u}}h)=q_{r,u}(h)$. Then for $\Psi \in \mathcal{H}^{2}(M,T)$, we have for every $r \geq 0$, $\omega \in \Omega$, $u \in U$:
\begin{eqnarray*}
\norm{\Psi(r,\omega,u) \circ i^{*}_{q_{r,u}} }^{2}_{\mathcal{L}_{2}(H_{q_{r,u}},G)} 
& = &  \sum_{k=1}^{\infty} q_{r,u}(i_{q_{r,u}} \circ \Psi(r,\omega,u)^{*}g_{k})^{2} \\
& = &  \sum_{k=1}^{\infty} q_{r,u}(\Psi(r,\omega,u)^{*}g_{k})^{2} \\
& = &  \sum_{k=1}^{\infty} \norm{Q^{1/2}_{r,u} \circ \Psi(r,\omega,u)^{*}g_{k}}^{2} \\
& = &  \norm{\Psi(r,\omega,u) \circ Q^{1/2}_{q_{r,u}}}^{2}_{\mathcal{L}_{2}(H,G)}.
\end{eqnarray*}
Now, Proposition \ref{propContEmbeddSquaIntegra} shows that 
for every $\Psi \in \mathcal{H}^{2}(M,T)$, we have $\Psi \circ i^{*}_{q} \in \Lambda^{2}(M,T)$. Hence \eqref{eqEquivaIntegNormCovarOper} follows from the equalities above and \eqref{eqDefiNormSquaIntegra}. 
\end{prf}

\begin{exam} \label{examStochIntegForLevyProcess}
Let $L=(L_{t}: t \geq 0)$ be a $H$-valued c\`{a}dl\`{a}g L\'{e}vy process. 
Let $U \in \mathcal{B}(H)$ be such that $0 \in U$ and $\int_{U} \, \norm{h}^{2} \nu(du)< \infty$ and let $\mathcal{R}= \{  U \cap \Gamma: \Gamma \in \mathcal{A} \} \cup \{ \{0\}\}$. Using the L\'{e}vy-It\^{o} decomposition of $L$ and the arguments from Example \ref{examLevyMartValuedMeasure}, we can define two independent martingale-valued measures $N_{1}$ and $N_{2}$ in $H$ by:
$$ N_{1}(t,A)=W_{t} \delta_{0}(A), \quad N_{2}(t,A)=\int_{A \setminus \{0 \}} \,  h \, \widetilde{N}(t,dh), \quad \forall \, t \geq 0, \, A \in \mathcal{R}. $$
with families of seminorms $p^{1}_{r,0}(h)^{2} = \inner{h}{Qh}$ and $p^{1}_{r,u}(h)^{2}=0$ if $u \neq 0$ for $N_{1}$, $p^{2}_{r,u}(h)^{2} = \inner{u}{h}^{2}$  for $N_{2}$; with the measures $\mu_{i} = \delta_{0}+ \restr{\nu}{U}$ and $\lambda_{i}=\mbox{Leb}$.

If $\Phi \in \Lambda^{2}(N_{1},T)$, we denote the stochastic integral $\displaystyle{\int_{0}^{t}\int_{U} \Phi(r,u)N_{1}(dr,du)}$ by $\displaystyle{\int_{0}^{t} \Phi(r,0) dW_{r}}$. Similarly, for $\Phi \in \Lambda^{2}(N_{2},T)$, we denote the stochastic integral $\displaystyle{\int_{0}^{t}\int_{U} \Phi(r,u)N_{2}(dr,du)}$ by $\displaystyle{\int_{0}^{t} \int_{U}  \Phi(r,u) \widetilde{N}(dt,du)}$. With the help of Proposition \ref{propCovarOperIntegrNorm} we can relate the above stochastic integrals with others theories available on the literature. 

In effect, suppose that $\Psi \in \mathcal{H}^{2}(N_{1},T)$.  
We know from Proposition \ref{propContEmbeddSquaIntegra} that $(\Psi(r,\omega,u) \circ i^{*}_{p^{1}_{r,u}}: r \in [0,T], \omega\in \Omega, u \in U ) \in \Lambda^{2}(N_{1},T)$ and hence the integral $\displaystyle{\int_{0}^{t} \Psi(r,0) dW_{r}}$ exists. On the other hand, since $\Psi \in \mathcal{H}^{2}(N_{1},T)$ the map  $(r,\omega) \mapsto \Psi(r,\omega,0)$ is $\mathcal{P}_{T}/\mathcal{B}(\mathcal{L}_{2}(H,G))$-measurable. Moreover, Proposition \ref{propCovarOperIntegrNorm} shows that  
$$\norm{\Psi}^{2}_{\Lambda^{2}(N_{1},T)}=\Exp \int_{0}^{T} \norm{ \Psi(r,0) \circ Q^{1/2}}^{2}_{\mathcal{L}_{2}(H,G)} dr < \infty.$$  
The above properties imply that the stochastic integral $\displaystyle{\int_{0}^{t} \Psi(r,0) dW_{r}}$ can be also defined using the theory of stochastic integration in Hilbert spaces as for example in  \cite{DaPratoZabczyk} or \cite{LiuRockner}. Since $W$ is a genuine Wiener process in $H$, our definition of the stochastic integral for simple processes coincides with the standard definition of stochastic integral in Hilbert spaces (see Remark \ref{remaCyliAndClassicIntegCoincide}). Therefore, both constructions provide the same stochastic process (up to indistinguishable versions). It is worth to mention that the integrability condition for the Wiener case can have another alternative expression. Indeed, if we let $\mathscr{H}= Q^{1/2}(H)$ equipped with the inner product $\inner{h_{1}}{h_{2}}_{\mathscr{H}} \defeq \inner{Q^{-1/2}h_{1}}{Q^{-1/2} h_{2}}$, then $\mathscr{H}$ is a Hilbert space and we have (see \cite{LiuRockner}, Sect. 2.3.2):
$$ \Exp \int_{0}^{T} \norm{ \Psi(r,0) \circ Q^{1/2}}^{2}_{\mathcal{L}_{2}(H,G)} dr= \Exp \int_{0}^{T} \norm{ \Psi(r,0)}^{2}_{\mathcal{L}_{2}(\mathscr{H},G)} dr. $$

Similarly, suppose that $\Psi \in \mathcal{H}^{2}(N_{2},T)$.  
From Proposition \ref{propContEmbeddSquaIntegra} we have $(\Psi(r,\omega,u) \circ i^{*}_{p^{2}_{r,u}}: r \in [0,T], \omega\in \Omega, u \in U ) \in \Lambda^{2}(N_{2},T)$ and hence the stochastic integral $\displaystyle{\int_{0}^{t}\int_{U} \Psi(r,u) \widetilde{N}(dr,du)}$ exists.
On the other hand, because $\Psi \in \mathcal{H}^{2}(N_{2},T)$ the  map $(r,\omega,u) \mapsto \Psi(r,\omega,u)$ is $\mathcal{P}_{T} \otimes \mathcal{B}(U) /\mathcal{B}(\mathcal{L}_{2}(H,G))$-measurable, and from Proposition \ref{propCovarOperIntegrNorm} we have that 
$$\norm{\Psi}^{2}_{\Lambda^{2}(N_{2},T)}=\Exp \int_{0}^{T} \int_{U} \norm{ \Psi(r,u) u }^{2} \nu(du) dr < \infty.$$  
Hence, the stochastic integral $\displaystyle{\int_{0}^{t}\int_{U} \Psi(r,u) \widetilde{N}(dr,du)}$ can be also defined using the theory of stochastic integration with respect to nuclear decomposable martingale-valued measures developed in \cite{Applebaum:2006}. The resulting integral coincide with ours since its action on simple integrands is the same (see Remark \ref{remaCyliAndClassicIntegCoincide}). 
\end{exam}

\begin{exam}\label{examItoIsomeCylLevySquaInteg} 
Let $L=(L_{t}: t \geq 0)$ be a cylindrical L\'{e}vy process in $H$ which has mean-zero and weak second moments.  
Let $M$ be the cylindrical martingale-valued measure defined by $L$ in Example \ref{examCylLevySquareMarting}. Since in this example we have $U=\{0\}$ and $\mu=\delta_{0}$, it is more convenient to denote the integrand $\Phi(r,\omega, 0)$ by $\Phi(r,\omega)$ and the integral $\displaystyle{\int_{0}^{t}\int_{U} \Phi(r,u) M(dr,du)}$ by $\displaystyle{\int_{0}^{t} \Phi(r) dL_{r}}$. 

With the help of Proposition \ref{propCovarOperIntegrNorm}, we can offer a detailed description of the integrands. To do this, let $\goth{Q}:H \rightarrow H$ be the non-negative, symmetric, continuous linear operator defined by 
$$ \inner{h_{1}}{\goth{Q}h_{2}}= \inner{h_{1}}{Qh_{2}}+ \int_{H} \inner{y}{h_{1}}\inner{y}{h_{2}} \nu(dy), \quad \forall h_{1}, h_{2} \in H. $$
A careful look at \eqref{eqCovariCylMartMeasLevySquar} shows that $\goth{Q}$ is indeed the cylindrical covariance of  $L$ and therefore we have $q_{r,0}(h_{1},h_{2})= \inner{h_{1}}{\goth{Q}h_{2}}$ for all $r \geq 0$. 

Then for every $\Psi \in \mathcal{H}^{2}(M,T)$ it follows from Proposition \ref{propCovarOperIntegrNorm} that 
\begin{eqnarray}
\norm{\Phi}^{2}_{\Lambda^{2}(M,T)} 
& = &\Exp \int_{0}^{T} \norm{ \Phi(r) \circ \goth{Q}^{1/2}}^{2}_{\mathcal{L}_{2}(H,G)} dr \nonumber \\
& = & \Exp \int_{0}^{T} \left[ \norm{ \Phi(r) \circ Q^{1/2}}^{2}_{\mathcal{L}_{2}(H,G)} + \int_{H } \norm{ \Phi(r) y }^{2} \nu(dy) \right] dr. \label{eqIntegrandsCylLevySquare} 
\end{eqnarray} 

Observe that since $\nu$ is only a cylindrical measure, the integral in the second line of \eqref{eqIntegrandsCylLevySquare} requires a clarification. A detailed explanation is given in page 14 of \cite{Riedle:2014}. 

Stochastic integration with respect to a  cylindrical L\'{e}vy process $L$ under the assumption of weak second moments was firstly considered in \cite{Riedle:2014}. In this work, the main class of integrands corresponds to the elements of the space $\mathcal{H}^{2}(M,T)$ (as defined in this example for $L$). An extension of the theory of integration in \cite{Riedle:2014}  to integrands having \eqref{eqIntegrandsCylLevySquare} finite is carried out in  \cite{KosmalaRiedle:VariaSPDE}. However, the reader must be warned that it might be possible that the stochastic integral constructed in this paper differs to those constructed in \cite{KosmalaRiedle:VariaSPDE, Riedle:2014}, being the main reason for this phenomena to occurs that the stochastic integral for the simple integrands in \cite{KosmalaRiedle:VariaSPDE, Riedle:2014} is carried out via a radonification of the increments which in principle is different form our approach based on Theorem \ref{theoRadonCylMartHilbSchmi} (see the discussion in Sect. \ref{sectIntro}). 
\end{exam}

\subsection{Further Properties of the Stochastic Integral}\label{secPropStochasInteg}

In this section we show that our definition of the stochastic integral satisfy all the standard properties. 

\begin{prop}\label{propMappingIntegContOpera} Let $E$ be a separable Hilbert space and let $R \in \mathcal{L}(G,E)$. Then, for each $\Phi \in \Lambda^{2}(M,T;H,G)$, we have $R \circ \Phi=\{ R \circ \Phi(r,\omega,u): r \in [0,T],  \omega \in \Omega, u \in U\} \in \Lambda^{2}(M,T;H,E)$, moreover $\Prob$-a.e. for every $t \in [0,T]$, 
\begin{equation} \label{eqIntUnderContMapping}
\int^{t}_{0} \int_{U} R \circ \Phi (r,u) M (dr, du)= R \left(\int^{t}_{0} \int_{U} \Phi (r,u) M (dr, du) \right). 
\end{equation} 
\end{prop}
\begin{prf} 
First, assume that $\Phi \in \mathcal{S}(M,T;H,G)$ is of the form \eqref{eqSimpleIntegrands}. Then, we have  
$$ R \circ \Phi(r,\omega,u)=  \sum_{i=1}^{n} \sum_{j=1}^{m}  \ind{]s_{i}, t_{i}]}{r} \ind{F_{i}}{\omega} \ind{A_{j}}{u} R \circ S_{i,j} \circ i^{*}_{q_{r,u}}, $$
and since each $R \circ S_{i,j} \in \mathcal{L}_{2}(H,E)$, it is clear that $R \circ \Phi \in 
 \mathcal{S}(M,T;H,E)$. Moreover, it follows from Definition \ref{defiClassSimpleInteg} that the stochastic integral takes the form 
$$ \int^{t}_{0} \int_{U} R \circ \Phi (r,u) M (dr, du)= \sum_{i=1}^{n} \sum_{j=1}^{m} \mathbbm{1}_{F_{i}} Y_{t}(s_{i},t_{i},A_{j},R \circ S_{i,j}).$$
However, observe from \eqref{eqDefRadonProcessY} that we have $\Prob$-a.e. for all $t \in [0,T]$ and $e \in E$,
$$ \inner{Y_{t}(s_{i},t_{i},A_{j},R \circ S_{i,j})}{e} = M((s_{i} \wedge t, t_{i} \wedge t],A_{j})(S_{i,j}^{*}\circ R^{*} e) = \inner{Y_{t}(s_{i},t_{i},A_{j},S_{i,j})}{R^{*}e}. $$
Therefore, $\Prob$-a.e. for all $t \in [0,T]$ and $e \in E$, we have
\begin{eqnarray*}
\inner{ \int^{t}_{0} \int_{U} R \circ \Phi (r,u) M (dr, du)}{e} 
& = & \inner{ \int^{t}_{0} \int_{U} \Phi (r,u) M (dr, du)}{R^{*}e} \\
& = & \inner{ R \left(\int^{t}_{0} \int_{U} \Phi (r,u) M (dr, du) \right)}{e}.
\end{eqnarray*}
But since $E$ is a separable Hilbert space, the above equality shows that \eqref{eqIntUnderContMapping} holds for $\Phi \in \mathcal{S}(M,T;H,G)$. 

Suppose now that $\Lambda^{2}(M,T;H,G)$. Then, we have 
\begin{enumerate}
\item $R \circ \Phi(r,\omega,u) \in \mathcal{L}_{2}(H_{q_{r,u}}, E)$, for all $r \in [0, T]$, $\omega \in \Omega$, $u \in U$, 
\item for each $h \in H$, $e \in E$, we have $R^* e \in G$ and hence the mapping 
$$(r,\omega,u) \mapsto q_{r,u}((R \circ \Phi(r,\omega,u))^* e, i_{q_{r,u}} h)= q_{r,u}(\Phi(r,\omega,u)^*( R^* e), i_{q_{r,u}} h),$$
is $\mathcal{P}_{T} \otimes \mathcal{B}(U)/\mathcal{B}(\R_{+})$-measurable.
\item \begin{multline*} 
\Exp \int_{0}^{T} \int_{U}  \norm{ R \circ \Phi(r,u)}_{\mathcal{L}_{2}(H_{q_{r,u}}, E)}^{2} \mu(du) \lambda(dr) \\
\leq \norm{R}_{\mathcal{L}(G, E)}^{2} \Exp \int_{0}^{T} \int_{U}  \norm{ \Phi(r,u)}_{\mathcal{L}_{2}(H_{q_{r,u}}, G)}^{2} \mu(du) \lambda(dr) < \infty.
\end{multline*} 
\end{enumerate}
Thus, we conclude that $R \circ \Phi \in \Lambda^{2}(M,T;H,E)$. Finally, since \eqref{eqIntUnderContMapping} holds for integrands in $\mathcal{S}(M,T;H,G)$, it follows from Proposition \ref{propSimpleDenseInSquareIntegra} and Theorem \ref{theoExisIntegralSquareInteg} that \eqref{eqIntUnderContMapping} also holds for integrands in $\Lambda^{2}(M,T;H,G)$. 
\end{prf}

\begin{prop} \label{propIntegralInSubintervalAndRandomSubset}
Let $0 \leq s_{0} < t_{0} \leq T$ and $F_{0} \in \mathcal{F}_{s_{0}}$. Then, for every $\Phi \in \Lambda^{2}(M,T)$, $\Prob$-a.e. we have  $\forall \, t \in [0,T]$, 
\begin{multline} \label{eqIntegralInSubintervalAndRandomSubset}
\int_{0}^{t} \int_{U} \, \mathbbm{1}_{]s_{0},t_{0}]\times F_{0}} \Phi(r,u) M(dr,du) \\
= \mathbbm{1}_{F_{0}} \left( \int_{0}^{t \wedge t_{0}} \int_{U} \, \Phi(r,u) M(dr,du) -\int_{0}^{t \wedge s_{0}}  \int_{U} \, \Phi(r,u) M(dr,du) \right). 
\end{multline}
\end{prop}
\begin{prf} First, by using similar ideas to those used in the proof of Proposition \ref{propMappingIntegContOpera} one can easily check that \eqref{eqIntegralInSubintervalAndRandomSubset} holds for $\Phi \in \mathcal{C}(M,T)$ of the  form \eqref{simpleFamiliesProofDenseSupspaceIntegrands}. Since $\mathcal{C}(M,T)$ spans $\mathcal{S}(M,T)$, the linearity of the integral shows that \eqref{eqIntegralInSubintervalAndRandomSubset} is valid for $\Phi \in \mathcal{S}(M,T)$. Finally, that \eqref{eqIntegralInSubintervalAndRandomSubset} holds for $\Phi \in \Lambda^{2}(M,T)$ can be show using the density of $\mathcal{S}(M,T)$  and the continuity of the stochastic integral. 
\end{prf}
 
\begin{prop} \label{propStoppedIntegral}
Let $\Phi \in \Lambda^{2}(M,T)$ and $\sigma$ be a $(\mathcal{F}_{t})$-stopping time such that $\Prob (\sigma \leq T)=1$. Then, $\Prob$-a.e. for every $t \in [0,T]$,
\begin{equation} \label{eqStoppedIntegral}
\int_{0}^{t} \int_{U} \, \ind{[0,\sigma]}{r} \Phi(r,u) M(dr,du) = \int_{0}^{t \wedge \sigma} \int_{U} \,  \Phi(r,u) M(dr,du).
\end{equation}
\end{prop}
\begin{prf} The proof can be carried out by following similar arguments to those in the proof of Lemma 2.3.9 in \cite{LiuRockner}, that is, by firstly checking \eqref{eqStoppedIntegral} for simple integrands of the form \eqref{eqSimpleIntegrands} and elementary stopping times. This can be done following similar arguments to those in the proof of Proposition \ref{propMappingIntegContOpera}. Then, using the density of $\mathcal{S}(M,T)$ and the continuity of the stochastic integral we can show that \eqref{eqStoppedIntegral} holds for $\Phi \in \Lambda^{2}(M,T)$. 
\end{prf}

The following result shows that the stochastic integral is linear on the integrators which are independent. Recall that for any two continuous Hilbertian seminorms $p$ and $q$ on $H$, with $p \leq q$, $i_{p,q}$ denotes the canonical inclusion from $H_{q}$ into $H_{p}$ (see Sect. \ref{sectionPrelim}).

\begin{prop} \label{propDecompSumIndepMartValMeasu}
Let $N_{1}$, $N_{2}$ be two cylindrical martingale-valued measures on $\R_{+} \times \mathcal{R}$, each satisfying Definition \ref{defiMartingValuedMeasure} determined by the family $\{ p_{r,u}^{j}: r \in \R_{+}, u \in U \}$ of continuous Hilbertian semi-norms on $H$ and measures $\lambda_{j}=\lambda$, $\mu_{j}=\mu$, for $j = 1,2$. Assume furthermore that for all $A, B \in \mathcal{R}$, the cylindrical processes $( N_{1}(t,A): t \geq 0)$ and $( N_{2}(t,B): t \geq 0)$ are independent. Let $M= (M(t,A): r \geq 0, A \in \mathcal{R})$ be given by the prescription:
\begin{equation} \label{eqSumIndepMartValuedMeasures}
M(t,A) \defeq N_{1}(t,A)+N_{2}(t,A), \quad \forall \, t \in \R_{+}, \, A \in \mathcal{R}. 
\end{equation} 
Then, $M$ is also a cylindrical martingale-valued measure on $\R_{+} \times \mathcal{R}$, satisfying Definition \ref{defiMartingValuedMeasure} for the measures $\lambda$, $\mu$ and the family of continuous Hilbertian semi-norms $\{ q_{r,u} : r \in \R_{+}, u \in U  \}$ satisfying $q_{r,u}(h)^{2}=p^{1}_{r,u}(h)^{2}+ p^{2}_{r,u}(h)^{2}$ for all $r \geq 0$, $u \in U$, $h \in H$. Moreover, if $\Phi \in \Lambda^{2}(M,T)$ we have:
\begin{enumerate}
\item For each $j=1,2$, $(\Phi(r,\omega,u) \circ i^{*}_{p_{r,u}^{j},q_{r,u}}: r \in [0,T], \omega \in \Omega, u \in U) \in \Lambda^{2}(N_{j},T)$, 
\item  $\Prob$-a.e., for all $t \in [0,T]$ we have, 
\begin{eqnarray*} 
\int_{0}^{t}\int_{U} \Phi(r,u) M(dr,du) 
& = &  \int_{0}^{t}\int_{U} \Phi(r,u) \circ i^{*}_{p_{r,u}^{1},q_{r,u}} N_{1}(dr,du) \\
& + & \int_{0}^{t}\int_{U} \Phi(r,u)\circ i^{*}_{p_{r,u}^{2},q_{r,u}} N_{2}(dr,du).  
\end{eqnarray*}
\end{enumerate}
\end{prop}

The proof of the above proposition follows from similar arguments to those used in the proof of Proposition 3.7 in \cite{BojdeckiJakubowski:1990} and for that reason we omit it. However, we illustrate its importance in the light of the following example. 

\begin{exam} \label{examSumLevyWienAndPoissoInteg}
Let $L=(L_{t}: t \geq 0)$ be a $H$-valued c\`{a}dl\`{a}g L\'{e}vy process and let $M$ be the L\'{e}vy martingale-valued measure defined in Example \ref{examLevyMartValuedMeasure}. Then, if $\Phi \in \Lambda^{2}(M,T)$ the stochastic integral 
$\displaystyle{\int_{0}^{t}\int_{U} \Phi(r,u) M(dr,du)}$ exist. By using Proposition \ref{propDecompSumIndepMartValMeasu}, we can relate this integral with those introduced in Example \ref{examStochIntegForLevyProcess}. 

In effect, let $N_{1}$ and $N_{2}$ be the two independent martingale-valued measures  in $H$ introduced in Example \ref{examStochIntegForLevyProcess}. It follows from the arguments given in Examples \ref{examLevyMartValuedMeasure} and \ref{examStochIntegForLevyProcess} that $M$, $N_{1}$ and $N_{2}$ satisfy \eqref{eqSumIndepMartValuedMeasures}. Moreover, it follows from \eqref{defiSemiNormsLevyMartValuedMeas} that $q_{r,0}=p^{1}_{r,0}$ and $q_{r,u}=p^{2}_{r,u}$ for $u \neq 0$. Therefore we have $H_{q_{r,0}}=H_{p^{1}_{r,0}}$ and $H_{q_{r,u}}=H_{p^{2}_{r,u}}$ for $u \neq 0$. Then for each $j=1,2$ the mapping $i^{*}_{p_{r,u}^{j},q_{r,u}} $ is nothing but the identity operator on $H_{q_{r,u}}$.  

Then if $\Phi \in \Lambda^{2}(M,T)$ it follows from Proposition \ref{propDecompSumIndepMartValMeasu} that $\Phi \in \Lambda^{2}(N_{j},T)$ for $j=1,2$, and that $\Prob$-a.e., for all $t \in [0,T]$ we have
\begin{equation}\label{eqLevyIntegAsSumWienerAndCompePoisson}
\int_{0}^{t}\int_{U} \Phi(r,u) M(dr,du)  =  \int_{0}^{t} \Phi(r,0)dW_{r} + \int_{0}^{t}\int_{U} \Phi(r,u) \widetilde{N}(dr,du),
\end{equation}
where the two integrals in the right-hand side of \eqref{eqLevyIntegAsSumWienerAndCompePoisson} correspond to the stochastic integrals with respect to $N_{1}$ and $N_{2}$ respectively (see Example  \ref{examStochIntegForLevyProcess}). Moreover, 
\begin{equation}
\norm{\Phi}^{2}_{\Lambda^{2}(M,T)}= \Exp \int_{0}^{T} \left[ \norm{ \Phi(r) \circ Q^{1/2}}^{2}_{\mathcal{L}_{2}(H,G)} + \int_{U} \norm{ \Phi(r) y }^{2} \nu(dy) \right] dr. \label{eqIntegraLevyMartiMeasure} 
\end{equation} 
\end{exam}

Our theory of stochastic integration also satisfies an stochastic Fubini theorem. In order to formulate our result we introduce the following class of integrands. 

Let $(E, \mathcal{E}, \varrho)$ be a $\sigma$-finite measure space. We denote by $\Xi^{1,2}(M,T,E)$ the linear space of all the families $\Phi=(\Phi(r,\omega,u,e): r \in [0,T], \omega \in \Omega, u \in U, e \in E)$ of operator-valued maps satisfying the following conditions:
\begin{enumerate}
\item $\Phi(r, \omega, u, e) \in \mathcal{L}_{2}(H_{q_{r,u}},G)$, $\forall r \in [0,T]$, $\omega \in \Omega$, $u \in U$, $e \in E$.   
\item The map $[0,T] \times \Omega \times U \times E \rightarrow \R_{+}$ given by $(r, \omega, u, e)\mapsto q_{r,u} (\Phi(r, \omega, u, e)^{*}g, i_{q_{r,u}} h)$ is $\mathcal{P}_{T} \otimes \mathcal{B}(U) \otimes \mathcal{E}/ \mathcal{B}(\R_{+})$-measurable, for every $h \in H$ and $g \in G$.
\item 
$$ \abs{\norm{\Phi}}_{\Xi^{1,2}(M,T,E)} \defeq \int_{E} \norm{\Phi(\cdot,\cdot,\cdot,e)}_{\Lambda^{2}(M,T)} \varrho (de) < \infty. $$
\end{enumerate}

For the above class of integrands we have the following: 

\begin{theo}[Stochastic Fubini's Theorem]\label{theoStochasticFubiniTheorem} Let $\Phi \in \Xi^{1,2}(M,T,E)$. Then,
\begin{enumerate}
\item For a.e. $(r,\omega, u) \in [0,T] \times \Omega \times U$, the mapping $E \backin e \mapsto \Phi(r,\omega,u,e) \in \mathcal{L}_{2}(H_{q_{r,u}},G)$ is Bochner integrable. Moreover,  
$$ \int_{E} \Phi(\cdot,\cdot,\cdot,e)\, \varrho \, (de)= \left\{ \int_{E} \Phi(r, \omega, u, e) \, \varrho (de): r \in [0,T], \omega \in \Omega,u \in U \right\} \in \Lambda^{2}(M,T;H,G) .$$ 
\item The mapping $E \backin e \mapsto  I(\Phi(\cdot,\cdot,\cdot,e)) =\left\{ \int_{0}^{t} \int_{U} \Phi(r,u,e) M(dr,du): t \in [0,T] \right\}  \in  \mathcal{M}^{2}_{T}(G)$ is Bochner integrable. Furthermore, $\Prob$-a.e. $\forall t \in [0,T]$, 
\begin{equation*} \label{defiIntegSecondPartStochFubini}
\left( \int_{E} \, I(\Phi(\cdot,\cdot,\cdot,e)) \, \varrho \, (de) \right)_{t} = \int_{E}  \left( \int_{0}^{t} \int_{U} \Phi(r,u,e) M(dr,du) \right) \varrho \, (de).
\end{equation*}
\item The following equality holds $\Prob$-a.e. $\forall t \in [0,T]$, 
\begin{equation} \label{identityFubiniTheoremWeakCase}
\int_{0}^{t} \int_{U}  \left( \int_{E} \Phi (\cdot,\cdot,\cdot,e) \, \varrho  (de)\right) M(dr,du) = \int_{E} \left( \int_{0}^{t} \int_{U} \Phi(r,u,e) M(dr,du) \right) \varrho (de).
\end{equation}
\end{enumerate}
\end{theo}

The main idea behind the proof of Theorem \ref{theoStochasticFubiniTheorem} is to first check that $\abs{\norm{\cdot}}_{\Xi^{1,2}(M,T,E)}$ defines a norm in $\Xi^{1,2}(M,T,E)$ that makes it a Banach space and that a dense subset for it is the collection $ S(M,T,E)$ of all families $ \Phi= (\Phi(r,\omega,u,e): r \in [0,T], \omega \in \Omega, u \in U, e \in E )$ of operator-valued mappings of the form:
\begin{equation} \label{defiSimpleIntegFubini}
 \Phi(r,\omega,u,e)= \sum^{p}_{l=1} \sum^{n}_{i=1} \sum^{m}_{j=1} \ind{]s_{i},t_{i}]}{r} \ind{F_{i}}{w} \ind{A_{j}}{u} \ind{D_{l}}{e} S_{i,j,l} \circ i^{*}_{q_{r,u}},
\end{equation}
where $m$, $n$, $p \in \N$, and for $l=1, \dots, p$, $i=1, \dots, n$, $j=1, \dots,m$, $0\leq s_{i} < t_{i} \leq T$, $F_{i} \in \mathcal{F}_{s_{i}}$, $A_{j} \in \mathcal{R}$, $D_{l} \in \mathcal{E}$ and $S_{i,j,l} \in \mathcal{L}_{2}(H,G)$. The proof of the denseness follows from a two-step process. First, if we denote by $\Xi^{2,2}(M,T,E)$ the subspace of all the $\Phi \in \Xi^{1,2}(M,T,E)$ satisfying
$\norm{\Phi}_{\Xi^{2,2}(M,T,E)} \defeq \int_{E} \norm{\Phi(\cdot,\cdot,\cdot,e)}^{2}_{\Lambda^{2}(M,T)} \varrho (de) < \infty$, then one can show that every element in $\Xi^{1,2}(M,T,E)$ can be approximated by a norm-increasing sequence of elements in $\Xi^{2,2}(M,T,E)$ (see the proof Lemma 4.22 in \cite{FonsecaMora:2018-1}). For the second step we must show that $ S(M,T,E)$ is $\norm{\cdot}_{\Xi^{2,2}(M,T,E)}$-dense in $\Xi^{2,2}(M,T,E)$, but this can be done following similar arguments to those used in the proof of Proposition \ref{propSimpleDenseInSquareIntegra}. 

Then for $\Phi \in  S(M,T,E)$ it is easy to verify that the conclusions of Theorem \ref{theoStochasticFubiniTheorem} are satisfied. The result extends to $\Phi \in \Xi^{1,2}(M,T,E)$ by following an approximation procedure with a sequence of elements in $S(M,T,E)$. A formal proof can be carried out by following the arguments used in the proof of Theorem 4.24 in \cite{FonsecaMora:2018-1}, but with the corresponding (mild) modifications to adapt it to our construction of the stochastic integral. For that reason we omit the proof and leave the details to the reader.

\subsection{An Extension of The Stochastic Integral}\label{sectExtenStocInteg}

As a final step in our construction of the stochastic integral, we extend the stochastic integral to the following class of integrands:

\begin{defi} \label{integrandsWeakIntegAlmostSureSquareMoments}
Let $\Lambda^{2,loc}(M,T)$ denote the collection of families $\Phi=(\Phi(r,\omega,u): r \in [0,T], \omega \in \Omega, u \in U )$ of operator-valued maps satisfying conditions (1) and (2) in Definition \ref{defiIntegrandsWithSquareMoments}, and satisfying: 
\begin{equation} \label{eqAlmostSureSecondMoment}
\Prob \left( \omega \in \Omega: \int_{0}^{T} \int_{U} \norm{ \Phi(r,u)}_{\mathcal{L}_{2}(H_{q_{r,u}}, G)}^{2} \mu(du) \lambda(dr) < \infty \right)=1.
\end{equation} 
\end{defi}

The space $\Lambda^{2,loc}(M,T)$ will be equipped with the vector topology $\mathcal{T}_{2,loc}^{M}$ generated by the local base of neighbourhoods of zero $\{ \Gamma_{\epsilon, \delta}: \epsilon > 0, \delta >0 \}$, given by 
$$ \Gamma_{\epsilon, \delta} = \left\{ \Phi \in \Lambda^{2,loc}(M,T): \Prob \left( \omega \in \Omega: \int_{0}^{T} \int_{U} \norm{ \Phi(r,u)}_{\mathcal{L}_{2}(H_{q_{r,u}}, G)}^{2} \mu(du) \lambda(dr) > \epsilon \right) \leq \delta \right\}. $$
Following similar arguments to those used in the proof of Proposition 2.4 in \cite{BojdeckiJakubowski:1990}, we can show that the space $(\Lambda^{2,loc}(M,T), \mathcal{T}_{2,loc}^{M})$ is a complete, metrizable topological vector space.  

Let $\Phi \in \Lambda^{2,loc}(M,T)$. If for each $n \in \N$ we define 
$$ \tau_{n} = \inf \left\{ t \in [0,T]: \int_{0}^{t} \int_{U} \norm{ \Phi(r,u)}_{\mathcal{L}_{2}(H_{q_{r,u}}, G)}^{2} \mu(du) \lambda(dr) \geq n \right\} \wedge T, $$
one can easily check that $( \tau_{n} : n \in \N)$ is an increasing sequence of $(\mathcal{F}_{t})$-stopping times satisfying \begin{inparaenum}[(i)] \item $\lim_{n \rightarrow \infty} \tau_{n} = T$ $\Prob$-a.e., and \item for each $n \in \N$, $\mathbbm{1}_{[0,\tau_{n}]} \Phi \in \Lambda^{2}(M,T)$. \end{inparaenum} A standard argument as in  Section 4.2 of \cite{DaPratoZabczyk} together with Theorem \ref{theoExisIntegralSquareInteg} and Proposition \ref{propStoppedIntegral} shows the existence of a unique  $G$-valued process 
$$\left( \int_{0}^{t} \int_{U} \, \Phi(r,u) M(dr,du) : t \in [0,T] \right)$$ 
that belongs to the space $\mathcal{M}^{2,loc}_{T}(G)$ of $G$-valued c\`{a}dl\`{a}g locally square integrable martingales, and such that for any increasing sequence of $(\mathcal{F}_{t})$-stopping times $( \sigma_{n} : n \in \N)$ satisfying $(i)$ and $(ii)$ given above, we have:
\begin{equation} \label{eqCompatibilityCondIntegral}
\int_{0}^{t \wedge \sigma_{n} }\int_{U} \, \Phi(r,u) M(dr,du) = \int_{0}^{t} \int_{U} \,  \mathbbm{1}_{[0,\sigma_{n}]} \Phi(r,u) M(dr,du), \quad \forall \, t \in [0,T], \, n \in \N.
\end{equation}
It should be clear that the property \eqref{eqCompatibilityCondIntegral} allow us to extend  to integrands in $\Lambda^{2,loc}(M,T)$ the properties of the stochastic integral listed in Section \ref{secPropStochasInteg}. Furthermore, if the space $\mathcal{M}^{2,loc}_{T}(G)$ is considered equipped with the topology of convergence in probability uniformly on $[0,T]$, then similar arguments to those in the proof of Proposition 4.20 in \cite{FonsecaMora:2018-1}
show that the map 
$$ \Lambda^{2,loc}(M,T) \ni \Phi \mapsto \left( \int_{0}^{t} \int_{U} \, \Phi(r,u) M(dr,du) : t \in [0,T] \right) \in \mathcal{M}^{2,loc}_{T}(G), $$
is linear and continuous.

\section{Stochastic Partial Differential Equations} \label{sectStocEvoluEquations}

\subsection{Existence and Uniqueness of Solutions}

In this section we apply our theory of stochastic integration to study stochastic partial differential equations of the form:
\begin{equation}\label{eqGeneralSEE}
\begin{cases}
dX_{t}  =(AX_{t}+B(t,X_{t}))dt+\int_{U} \, F(t,u,X_{t}) M(dt,du), \\
X_{0}  = \xi. 
\end{cases}
\end{equation}
Here we assume the following:
\begin{assu}\label{assuCoefficiSEE} \hfill

\begin{enumerate}
\item $\xi$ is a $\mathcal{F}_{0}$-measurable $G$-valued random variable. 
\item $A: \mbox{Dom}(A) \subseteq G \rightarrow G$ is the infinitesimal generator of a $C_{0}$-semigroup $(S(t): t \geq 0)$ on $G$. 
\item $B: [0,T] \times \Omega \times G \rightarrow G$ is $\mathcal{P}_{T} \otimes \mathcal{B}(G)/ \mathcal{B}(G)$ measurable. 
\item $F=(F(r,u,g): r \in \R_{+}, u \in G, g \in G)$ is such that:
\begin{enumerate}
\item $F(r,y,g) \in \mathcal{L}_{2}(H_{q_{r,u}},G)$ $\forall r \geq 0$, $u \in U$, $g \in G$, 
\item The mapping $(r,u, g) \mapsto q_{r,u}(F(r,u,g_{1})^* g_{2}, i_{q_{r,u}} h)$ is $\mathcal{B}(\R_{+}) \otimes \mathcal{B}(U) \otimes \mathcal{B}(G)/\mathcal{B}(\R_{+})$-measurable, $\forall g_{1}, g_{2} \in G$, $h \in H$. 
\end{enumerate}
\item $M$ is a cylindrical martingale-valued measure is as in Definition \ref{defiMartingValuedMeasure} with the particularity that $\lambda$ is the Lebesgue measure on $\R_{+}$. 
\end{enumerate}
\end{assu}

Recall that since $G$ is a Hilbert space, the collection $(S(t)^{*}: t \geq 0)$ defines also a $C_{0}$-semigroup on $G$ whose infinitesimal generator is given by $A^{*}$ (see Corollary 1.10.6 in \cite{Pazy}, p.41). 

A predictable $G$-valued process $(X_{t}: t \geq 0)$ is called a \emph{mild solution} to \eqref{eqGeneralSEE} if $\forall t \geq 0$, $\Prob$-a.e.
\begin{equation} \label{eqMildSolGeneralSEE}
X_{t}=S(t)\xi +\int_{0}^{t} S(t-r)B(r, X_{r})dr+\int_{0}^{t} \int_{U} S(t-r) F(r,u,X_{r})M(dr,du),
\end{equation}
where in the right-hand side in \eqref{eqMildSolGeneralSEE} the first integral is defined pathwise as a Bochner integral and the second integral as a stochastic integral, provided both of them are well-defined. 

We will assume the following growth and Lipschitz conditions in our coefficients $B$ and $F$:

\begin{assu}\label{assuGrowthLipscCondi} 
 There exists two functions $a, b: \R_{+} \rightarrow \R_{+}$ satisfying $\displaystyle{ \int_{0}^{T} a(r)+ b(r)^{2} dr < \infty}$ for each $T>0$, and such that: 
\begin{enumerate}
\item (Growth conditions) For all $r>0$, $g \in G$:
\begin{eqnarray*}
\norm{B(r, g)} & \leq &  a(r) (1+\norm{g}), \\
\int_{U} \, \norm{F(r, u, g)}^{2}_{\mathcal{L}_{2}(H_{q_{r,u}}, G) } \mu(du) & \leq &  b(r)^{2} (1+\norm{g})^{2}.  
\end{eqnarray*}
\item (Lipschitz conditions) For all $r>0$, $g_{1}$, $g_{2} \in G$:
\begin{eqnarray*}
\norm{B(r, g_{1})- B(r, g_{2})} & \leq &  a(r) \norm{g_{1}-g_{2}}, \\
\int_{U} \, \norm{F(r, u, g_{1})- F(r, u, g_{2})}^{2}_{\mathcal{L}_{2}(H_{q_{r,u}}, G) } \mu(du)  & \leq &  b(r)^{2} \norm{g_{1}-g_{2}}^{2}. 
\end{eqnarray*}
\end{enumerate}
\end{assu}

\begin{rema}\label{remaStandCondCoefficBandF}
Suppose that the coefficient $F$ is of the form $F:\R_{+} \times U \times G \rightarrow \mathcal{L}_{2}(H,G)$.  
In such a case and by using Proposition \ref{propContEmbeddSquaIntegra} the reader can check that the growth and Lipschitz conditions in Assumption \ref{assuGrowthLipscCondi}(2) are implied by the following standard growth and Lipschitz conditions: 
\begin{eqnarray*}
 \int_{U} \norm{F(r, u, g)}_{\mathcal{L}_{2}(H, G) } \mu(du) & \leq &  b(r) (1+\norm{g}), \\
 \int_{U} \norm{F(r, u, g_{1})- F(r, u, g_{2})}_{\mathcal{L}_{2}(H, G) } \mu(du)  & \leq &  b(r) \norm{g_{1}-g_{2}}.  
\end{eqnarray*}
\end{rema}

The following result is the main existence and uniqueness result: 

\begin{theo}\label{theoExisUniqSoluGenerSEE}
Assume that $\xi$ is square integrable. Then, there exist a unique (up to modifications) $G$-valued predictable process $X=(X_{t}: t \geq 0 )$ that is a mild solution to \eqref{eqGeneralSEE} satisfying $\displaystyle{\sup_{t\in[0,T]}\Exp \norm{X_{t}}^2<\infty}$ for each $T>0$.
\end{theo}
\begin{prf} It is sufficient to show that the result holds on the bounded interval $[0,T]$ for any given $T>0$. 
Let $\Upsilon_{T}$ denotes the Banach space of all (equivalence classes) of $G$-valued predictable processes $X=(X_{t}: t \in [0,T])$ such that 
$$ \norm{X}_{T} \defeq \left( \sup_{t \in [0,T]} \Exp \norm{X_{t}}^{2} \right)^{1/2} < \infty.$$
Define the operator 
$K:\Upsilon_{T} \rightarrow \Upsilon_{T}$ as $K(X)=K_{0}(X)+K_{1}(X)+K_{2}(X)$, where
\begin{align*}
& K_{0}(X)=S(t)\xi \\
& K_{1}(X)=\int_{0}^{t} S(t-r)B(r, X_{r})dr \\
& K_{2}(X)=\int_{0}^{t} \int_{U} S(t-r) F(r,u,X_{r})M(dr,du). 
\end{align*}
We will check that $K$ is well-defined. To do this, let $X \in \Upsilon_{T}$. First, recall that since $S(t)$ is a $C_{0}$-semigroup, there exists $\alpha>0$, $N \geq 1$ such that $\norm{S(t)g} \leq N e^{\alpha t} \norm{g} $, $\forall t \geq 0$, $g \in G$ (see e.g. \cite{Pazy}). The strong continuity of the $C_{0}$-semigroup shows that $K_{0}(X)$ is $(\mathcal{F}_{t})$-adapted and continuous, hence has a predictable version. Moreover, 
$$ \sup_{t \in [0,T]} \Exp \norm{K_{0}(X)_{t}}^{2} \leq N^{2} e^{2\alpha T} \Exp \norm{\xi}^{2}< \infty. $$
For $K_{1}(X)$, observe that from the growth condition in Assumption \ref{assuGrowthLipscCondi}(1), it is easy to verify that 
\begin{equation} \label{eqSquaMomenK1}
\sup_{t \in [0,T]} \Exp \norm{K_{1}(X)}^{2} \leq 2 N^{2} e^{2\alpha T} \left( \int_{0}^{T} a(r) dr \right)^{2} \left( 1+ \sup_{t \in [0,T]} \Exp \norm{X_{t}}^{2} \right) < \infty.  
\end{equation} 
Then, the Bochner integral in $K_{1}(X)$ is $\Prob$-a.e. well defined for each $t \in [0,T]$. Moreover, similar arguments to those in the proof of Lemma 6.2.9 in \cite{LiuRockner} shows that this process is continuous. Since it is clearly $(\mathcal{F}_{t})$-adapted, it has a predictable version. Then, by \eqref{eqSquaMomenK1} we have $K_{1}(X) \in \Upsilon_{T}$. 

For $K_{2}(X)$, we need to show that for any given $t \in [0,T]$ we have 
$$\displaystyle{\left( \ind{[0,t]}{r} S(t-r) F(r,u,X_{r}(\omega)): r \in [0,T], \omega \in \Omega, u \in U \right) \in \Lambda^{2}(M,t)}.$$ 
In effect, it is clear that $S(t-r) F(r,u,X_{r}(\omega))  \in \mathcal{L}_{2}(H_{q_{r,u}},G)$ for each $r \in [0,T]$, $\omega \in \Omega$, $u \in U$. Now, from the strong continuity of the dual semigroup $(S(t)^{*}: t \geq 0)$, the predictability of $X$ and the $q_{r,u}$-measurability of $F$, we have that the map
$$ (r,\omega,u) \mapsto q_{r,u}(F(r,u,X_{r}(\omega))^* S(t-r)^{*} g, i_{q_{r,u}} h),$$ 
is $\mathcal{P}_{T} \otimes \mathcal{B}(U)/\mathcal{B}(\R_{+})$-measurable for each $h \in H$, $g \in G$. Furthermore, from the growth condition Assumption \ref{assuGrowthLipscCondi}(2) we have
\begin{flalign}
 & \Exp \int_{0}^{t} \int_{U}  \norm{ S(t-r) F(r,u, X_{r})}_{\mathcal{L}_{2}(H_{q_{r,u}}, G)}^{2} \mu(du) dr \label{eqSeconMomeStochConv} \\
 & \leq 2 N^{2} e^{2 \alpha t} \left( \int_{0}^{T} b(r)^{2} dr \right) \left( 1+ \sup_{t \in [0,T]} \Exp \norm{X_{t}}^{2} \right)< \infty.  \nonumber
\end{flalign}
Hence, the stochastic integral defining $K_{2}(X)$ is well-defined. It is by construction $(\mathcal{F}_{t})$-adapted and using similar arguments to those in the proof of Theorem 6.14 in \cite{FonsecaMora:2018-1} we can show it is mean-square continuous, hence it has a predictable version (see \cite{DaPratoZabczyk}, Proposition 3.6, p.76). Furthermore, by the It\^{o} isometry and \eqref{eqSeconMomeStochConv} we have $K_{2}(X) \in \Upsilon_{T}$. Therefore, the operator $K$ is well-defined. 

Now, it is simple to check that 
$$ \norm{X}_{T, \beta} \defeq \left( \sup_{t \in [0,T]} e^{-\beta t} \Exp \norm{X_{t}}^{2} \right)^{1/2},$$
defines an equivalent norm on $\Upsilon_{T}$ for each $\beta>0$. Hence, following the proof of Theorem 9.29 in \cite{PeszatZabczyk}, to show existence and uniqueness of a mild solution to \eqref{eqGeneralSEE}, it is enough to show that $K$ is a contraction under the norm $\norm{\cdot}_{T, \beta}$ for a suitable $\beta$. Then, the result will be a direct consequence of Banach's fixed point theorem. 

In effect, let $X, Y \in \Upsilon_{T}$. Then, by the Lipschitz condition in Assumption \ref{assuGrowthLipscCondi}(1) and following similar calculations to those in the proof of Theorem 9.29 in \cite{PeszatZabczyk}, we have
$$ \norm{K_{1}(X)-K_{1}(Y)}^{2}_{T, \beta} \leq   N^{2} e^{2\alpha T} \left( \int_{0}^{T} a(r) dr \right) \left( \int_{0}^{T} a(r) e^{-\beta r} dr \right) \norm{X-Y}^{2}_{T, \beta}. $$
On the other hand, by the Lipschitz condition in Assumption \ref{assuGrowthLipscCondi}(2), the It\^{o} isometry,  and again  following similar calculations to those in the proof of Theorem 9.29 in \cite{PeszatZabczyk}, we have
$$ \norm{K_{2}(X)-K_{2}(Y)}^{2}_{T, \beta} \leq   N^{2} e^{2\alpha T} \left( \int_{0}^{T} b(r)^{2} e^{-\beta r} dr \right) \norm{X-Y}^{2}_{T, \beta}. $$ 
Hence, from the two estimates above, we have 
$$ \norm{K(X)-K(Y)}_{T, \beta} \leq C_{T,\beta} \norm{X-Y}_{T, \beta}, $$
with
$$C_{T,\beta}^{2} =  2 N^{2} e^{2\alpha T} \left[ \left( \int_{0}^{T} b(r)^{2} e^{-\beta r} dr \right) + \int_{0}^{T} b(r)^{2} e^{-\beta r} dr  \right].$$
Taking $\beta$ sufficiently large such that $ C_{T,\beta}<1$, we obtain that $K$ is a $\norm{\cdot}_{T, \beta}$-contraction. Hence, by the arguments given above this shows the existence and uniqueness of a mild solution to \eqref{eqGeneralSEE} satisfying $\sup_{t\in[0,T]}\Exp \norm{X_{t}}^2<\infty$.
\end{prf}

We finalize this section by briefly discussing the existence of a weak solution to \eqref{eqGeneralSEE}.  

A predictable $G$-valued process $(X_{t}: t \geq 0)$ is called a \emph{weak solution} to \eqref{eqGeneralSEE} if $\forall t \geq 0$ and $g \in \mbox{Dom}(A^{*})$, $\Prob$-a.e.
\begin{equation} \label{eqWeakSolGeneralSEE}
\inner{X_{t}}{g}=\inner{\xi}{g} +\int_{0}^{t} \inner{X_{r}}{A^{*}g}+ \inner{B(r, X_{r})}{g} dr+\int_{0}^{t} \int_{U} \inner{F(r,u,X_{r})}{g} M(dr,du),
\end{equation}
where in the right-hand side in \eqref{eqWeakSolGeneralSEE} the first integral is defined pathwise as a Lebesgue integral and the second integral as a stochastic integral, provided both of them are well-defined. 

In the following result we provide sufficient conditions for the equivalence of weak and mild solutions. 

\begin{theo}\label{theoEquivWeakMildSolu}
Let $X=(X_{t}: t \geq 0 )$  be a $G$-valued predictable process and assume that for each $T>0$, we have for $X$, $B$ and $F$ that:
\begin{enumerate}
\item $\displaystyle{ \Exp \int_{0}^{T} \, \norm{X_{r}} \, dr < \infty.} $
\item $\displaystyle{ \Exp \int_{0}^{T} \, \norm{B(r,X_{r})} \, dr < \infty.} $
\item $\displaystyle{ \Exp \int_{0}^{T} \int_{U} \, \norm{F(r,u,X_{r})}^{2}_{\mathcal{L}_{2}(H_{q_{r,u}},G)} \mu(du) dr < \infty.} $
\end{enumerate} 
Then, $X$ is a weak solution to \eqref{eqGeneralSEE} if and only if $X$ is a mild solution to \eqref{eqGeneralSEE}. 
\end{theo}

A proof for Theorem \ref{theoEquivWeakMildSolu} can be carried out following very closely the arguments used in the proof of Theorem 6.9 in \cite{FonsecaMora:2018-1} (there in the context of the dual of a nuclear space) by using the stochastic Fubini theorem (Theorem \ref{theoStochasticFubiniTheorem}). For these reasons we decided not to include a proof here. 

Now, observe that if $X=(X_{t}: t \geq 0 )$ is the unique mild solution of Theorem \ref{theoExisUniqSoluGenerSEE} to \eqref{eqGeneralSEE}, then since $\sup_{t\in[0,T]}\Exp \norm{X_{t}}^2<\infty$ and the growth conditions on $B$ and $F$ imply the conditions in Theorem \ref{theoEquivWeakMildSolu}, then it follows that $X$ is also a weak solution to \eqref{eqGeneralSEE}. We state it formally for further references:

\begin{coro}\label{coroWeakSolu}
The unique mild solution $X=(X_{t}: t \geq 0 )$ of Theorem \ref{theoExisUniqSoluGenerSEE} to \eqref{eqGeneralSEE}
is also a weak solution to \eqref{eqGeneralSEE}. 
\end{coro}

\subsection{Stochastic Partial Differential Equations Driven by  L\'{e}vy Noise} \label{sectSEELevyNoise}

In this section we show that our theory of stochastic integration and Theorem \ref{theoExisUniqSoluGenerSEE} can be used to show existence and uniqueness of solutions to some classes of stochastic partial differential equations driven by a $H$-valued L\'{e}vy process $L=(L_{t}: t \geq 0)$ where the coefficients depends not only on the solution but also depends on the time and the jump-space variables. It is worth to stress that no finite moments assumptions on $L$ will be required. 

To set up our problem, assume that $L$ has the L\'{e}vy-It\^{o} decomposition \eqref{eqLevyItoDecomp}. In such a case, we can formulate a general stochastic partial differential equation driven by $L$ as follows:
\begin{multline}\label{eqLevyDrivenSEELID}
 d X_{t}= (A X_{t}+ B (t,X_{t})) dt+ F(t,0,X_{t}) dW_{t} \\
+ \int_{\norm{h}<1} F(t,h,X_{t}) \tilde{N}(dt,dh)+ \int_{\norm{h} \geq 1} F(t,h,X_{t}) N(dt,dh), 
\end{multline}
for all $t \geq 0$ with initial condition $X_{0}=\xi$. We assume $\xi$, $A$, $(S(t): t \geq 0)$, $B$ and $F$ satisfy Assumptions \ref{assuCoefficiSEE} (1)-(4) for $U= H$, $\mu=\nu$ where $\nu$ is the L\'{e}vy measure of $L$, $\lambda$ is the Lebesgue measure on $\R_{+}$, and with the family of continuous Hilbertian semi-norms  $( q_{r,u}: r \in \R_{+}, u \in H )$ defined in \eqref{defiSemiNormsLevyMartValuedMeas}.  Furthermore, we assume that $B$ and $F$ satisfy Assumption \ref{assuGrowthLipscCondi}. Observe that from the particular definition of the seminorms $( q_{r,u})$, the growth condition for $F$ now takes the form:
\begin{equation} \label{eqGrowthCondLevySEE}
\norm{ F(r,\omega, 0,g) \circ Q^{1/2}}^{2}_{\mathcal{L}_{2}(H,G)} + \int_{H} \norm{ F(r,\omega,h,g) h }^{2} \nu(du) \leq b(r)^{2}(1+\norm{g})^{2}, 
\end{equation}
and the Lipschitz condition is given by: 
\begin{multline} \label{eqLipschitzCondLevySEE}
\norm{ \left(F(r,\omega, 0,g_{1})-F(r,\omega, 0,g_{2}) \right) \circ Q^{1/2}}^{2}_{\mathcal{L}_{2}(H,G)} \\
+ \int_{H} \norm{ \left(F(r,\omega,h,g_{1}) -F(r,\omega,h,g_{2}) \right)h }^{2} \nu(du) \leq b(r)^{2}\norm{g_{1}-g_{2}}^{2}.  
\end{multline}
Now, we will define a solution to \eqref{eqLevyDrivenSEELID} by introducing a sequence of stopping times that transforms the equation \eqref{eqLevyDrivenSEELID} into a system of equations of the form \eqref{eqGeneralSEE} which can solved using Theorem \ref{theoExisUniqSoluGenerSEE}. This procedure is a modification of the one used by Peszat and Zabczyk in Section 9.7 in \cite{PeszatZabczyk}. However, our procedure will differ in that it uses the family of L\'{e}vy martingale-valued measures in Example \ref{examLevyMartValuedMeasure} and the linearity of the stochastic integral on its integrators given in Proposition \ref{propDecompSumIndepMartValMeasu}. 

For each $n \in \N$, let $U_{n} = \{ h \in H: \norm{h} < n \}$. Following Section 9.7 in \cite{PeszatZabczyk}, define a sequence of stopping times by   
\begin{equation} \label{defiStoppingTimesTauNLevyNoise}
 \tau_{n}(\omega) \defeq \inf \{ t \geq 0: \Delta L_{t} (\omega) \notin U_{n} \}, \quad \forall \, \omega \in \Omega.
\end{equation}  
Since $H = \bigcup_{n \in \N} U_{n}$, it follows that $\tau_{n} \rightarrow \infty$ as $n \rightarrow \infty$.

Let $\mathcal{R}= \mathcal{A} \cup \{0\}$, where $\mathcal{A}$ denotes the collection of all $A \in \mathcal{B}(H \setminus \{ 0 \})$ such that $0 \notin \overline{A}$. For every $n \in \N$, let $M_{n}=(M_{n}(t,A): r \geq 0, A \in \mathcal{R})$ be the L\'{e}vy martingale-valued measure given by
\begin{equation} \label{levyMartValuedMeasUn}
M_{n}(t,A) = W_{t} \delta_{0}(A) + \int_{U_{n} \cap (A \backslash \{0 \})} h \widetilde{N}(t,dh), \quad \mbox{ for } \, t \geq 0, \, A \in \mathcal{R}. 
\end{equation}
Observe that $M_{n}$ is well-defined (see Example \ref{examLevyMartValuedMeasure}) since 
\begin{equation*} \label{intLevyMeasuOnUnIsBounded}
\int_{U_{n}} \, \norm{h}^{2} \nu(du) \leq  \int_{\norm{h}< 1 } \norm{h}^{2} \nu(dh)+ n^{2}  \nu ( \{ h: \norm{h} \geq 1 \} ) < \infty. 
\end{equation*}
Define $B_{n}$ by 
$$B_{n}(t,\omega,g)=B(t,\omega,g) + \int_{U_{n}\setminus U_{1}} F(r,\omega,h,g) h \, \nu(dh).$$
The integral with respect to $\nu$ is well-defined as can be seen directly from \eqref{eqGrowthCondLevySEE}. Moreover, one can easily check that $B_{n}$ also satisfies the growth and Lipschitz conditions. Therefore, Theorem \ref{theoExisUniqSoluGenerSEE} and Corollary \ref{coroWeakSolu} show that for every $n \in \N$ the following abstract Cauchy problem
\begin{equation} \label{eqLevyAbsCauchyProbUn}
\begin{cases}
d X^{(n)}_{t}= (AX^{(n)}_{t}+ B_{n}(t,X^{(n)}_{t})) dt+\int_{U_{n}} F(t,u,X^{(n)}_{t}) M_{n} (dt,du), \quad \mbox{for }t \geq 0, \\
X^{(n)}_{0}=\xi.
\end{cases}
\end{equation}
has a unique predictable solution $X^{(n)}$ such that for each $T>0$, $\sup_{t \in [0,T]} \Exp \norm{X^{(n)}_{t}}^{2} < \infty $.

\begin{theo} \label{theoCompatLevyAbsCauchyProbUn}
For every $t \in [0,T]$ and all $m \leq n$, $X^{(n)}_{t}=X^{(m)}_{t}$ $\Prob$-a.e. on $\{ t \leq \tau_{m}\}$. Moreover, the $H$-valued predictable process $X$ defined by $X_{t}=X^{(m)}_{t}$ for $t \leq \tau_{m}$ is a mild and a weak solution to \eqref{eqLevyDrivenSEELID}. 
\end{theo}
\begin{prf} Our proof follows closely the proof of Theorem 7.1 in  \cite{FonsecaMora:2018-1}. For that reason we only provide the main arguments. 

Let $t \in [0,T]$. For $m \leq n$, the properties of the Poisson integrals shows that $M_{n}-M_{m}$ also a (cylindrical) martingale-valued measure on $\R_{+} \times \mathcal{R}$ and that is independent of $M_{n}$. Hence, because $X^{(n)}$ and $X^{(m)}$ are mild solutions to \eqref{eqLevyAbsCauchyProbUn} and from Proposition \ref{propDecompSumIndepMartValMeasu}, we have $\Prob$-a.e. 
\begin{eqnarray*}
X^{(n)}_{t} - X^{(m)}_{t}
& = & \int_{0}^{t} S(t-r)(B(r,X^{(n)}_{r})-B(r,X^{(m)}_{r})) dr \\
& {} & + \int_{0}^{t} \int_{U_{n}\setminus U_{1}} S(t-r) F(r,h,X^{(n)}_{r})h \, \nu (dh) dr \\
& {} & - \int_{0}^{t} \int_{U_{m}\setminus U_{1}} S(t-r) F(r,h,X^{(m)}_{r})h \, \nu (dh) dr \\
& {} & + \int_{0}^{t} \int_{H} S(t-r)(F(r,h,X^{(n)}_{r})-F(r,h,X^{(m)}_{r})) \, M_{m}(dr,dh) \\
& {} & + \int_{0}^{t} \int_{H} S(t-r) F(r,h,X^{(n)}_{r}) \, (M_{n}-M_{m})(dr,dh) 
\end{eqnarray*}
Now, it follows from \eqref{levyMartValuedMeasUn} that on the set $\{ t \leq \tau_{m}\}$ we have
\begin{equation*} \label{eqDescripMnMinusMn}
M_{n}(r,A)(f)-M_{m}(r,A)(f)
=-r \int_{(U_{n} \setminus U_{1}) \cap (A \backslash \{0 \})} \inner{f}{h} \nu(dh) + r \int_{(U_{m} \setminus U_{1}) \cap (A \backslash \{0 \})} \inner{f}{h} \nu(dh). 
\end{equation*}
Therefore,  
\begin{eqnarray*}
(X^{(n)}_{t} - X^{(m)}_{t}) \mathbbm{1}_{\{t \leq \tau_{m}\}} 
& =  & \int_{0}^{t} S(t-r)(B(r,X^{(n)}_{r})-B(r,X^{(m)}_{r})) dr \mathbbm{1}_{\{t \leq \tau_{m}\}} \\
& + & \int_{0}^{t} \int_{U_{m}\setminus U_{1}} S(t-r)(F(r,h,X^{(n)}_{r})-F(r,h,X^{(m)}_{r}))h  \, \nu (dh) dr \mathbbm{1}_{\{t \leq \tau_{m}\}} \\
& + & \int_{0}^{t} \int_{H} S(t-r)(F(r,h,X^{(n)}_{r})-F(r,h,X^{(m)}_{r})) \, M_{m}(dr,dh) \mathbbm{1}_{\{t \leq \tau_{m}\}}.
\end{eqnarray*}
Then, \eqref{defiSemiNormsLevyMartValuedMeas},  It\^{o} isometry \eqref{eqItoIsometrySquareIntegrands} and the Lipschitz conditions on $B$ and $F$ yields:
\begin{equation*}
Y(t)  \leq 3\, N^{2} e^{2 \alpha t} \int_{0}^{t} (a(r)^{2} +2b(r)^{2}) Y(r) dr, 
\end{equation*}
where $Y(t)= \Exp \left( \norm{X^{(n)}_{t} - X^{(m)}_{t}}^{2} \mathbbm{1}_{\{t \leq \tau_{m}\}} \right) $. Observe that $\sup_{t \in [0,T]} Y(t)< \infty$. Hence, from Gronwall's inequality if follows that $Y(t)=0$. Then, we conclude that $X^{(n)}_{t}=X^{(m)}_{t}$ $\Prob$-a.e. on $\{ t \leq \tau_{m}\}$. Finally, because $X^{(n)}$ is by definition a mild and a weak solution to \eqref{eqLevyDrivenSEELID} on $\{ t \leq \tau_{m}\}$, then $X$ is also a mild and a weak solution to \eqref{eqLevyDrivenSEELID}.    
\end{prf}
   
\begin{rema}
In \cite{PeszatZabczyk}, Section 9.7, the authors consider stochastic partial differential equations of the form
$$ dX_{t}=(AX_{t}+B(X_{t}))dt+F(X_{t})dM_{t}+G(X_{t})dP_{t},$$ 
where $M$ and $P$ are respectively the martingale-part and the compound Poisson process part of the L\'{e}vy-It\^{o} decomposition of a L\'{e}vy process in $H$. There, the authors were able to prove existence and uniqueness for coefficients which does not necessarily take values
in the space of Hilbert-Schmidt operators by utilizing growth and Lipschitz conditions that consider the smoothing property of the $C_{0}$-semigroup. Observe that in our growth and Lipschitz conditions \eqref{eqGrowthCondLevySEE} and \eqref{eqLipschitzCondLevySEE} we are not considering such a smoothing effect. This difference arises because in  \eqref{eqLevyDrivenSEELID} our coefficients depends also on the time and jump-space variable. In that sense, we believe that the result in our Theorem \ref{theoCompatLevyAbsCauchyProbUn} can serve as a complement to the theory developed in \cite{PeszatZabczyk}.  
\end{rema}

\smallskip

\textbf{Acknowledgements} { The authors acknowledge The University of Costa Rica for providing financial support through the grant ``B9131-An\'{a}lisis Estoc\'{a}stico con Procesos Cil\'{i}ndricos''. 

\end{document}